\documentclass[11pt,letter,subeqn,fleqn]{article}

\usepackage{fullpage}
\usepackage{hyperref}
\usepackage{amsmath, amssymb,amsthm,cite}
\usepackage{graphics,color}
\usepackage{graphicx}
\usepackage{tabulary}
\usepackage{float}
\usepackage{graphicx}
\usepackage{multirow}
\usepackage{subcaption}
\usepackage{longtable}

\begin{document}

\tolerance=9999

\def \dfracskip{2ex}
\def\beq{\begin{equation}}
\def\eeq{\end{equation}}
\def\barr{\begin{array}{ll}}
\def\earr{\end{array}}
\def\const{\hbox{\rm const}}

\def\sg#1{{\rm #1}}

\let\vector\vec %

\def\vec#1{{\boldsymbol{\rm #1}}}  %
\def\div{\mathop{\hbox{\rm div}}}
\def\grad{\mathop{\hbox{\rm grad}}}
\def\curl{\mathop{\hbox{\rm curl}}}

\def \ConsLawDensity{\Psi}
\def \ConsLawFlux{\Phi}

\def\depvars{\mathbf{V}}
\def\depvar{\depvars}

\def\Dop{\mathcal{D}}

\def\lieder#1{\mathcal{L}_{#1}}

\def\const{\hbox{\rm const}}
\def\diag{\hbox{\rm diag}}
\def\mod{\hbox{\rm mod}}
\def\rank{\mathop{\hbox{\rm rank}}}
\def\Image{\mathop{\hbox{\rm Image}}}
\def\mes{\mathop{\hbox{\rm mes}}}
\def\grad{\mathop{\hbox{\rm grad}}}
\def\sign{\mathop{\hbox{\rm sign}}}
\def\ind{\mathop{\hbox{\rm ind}}}
\def\id{\mathop{\hbox{\rm id}}}
\def\d{\mathop{\hbox{\rm d}}}
\def\ad{\hbox{ad}}
\def\ch{\mathop{\hbox{\rm ch}}}
\def\th{\mathop{\hbox{\rm th}}}
\def\cth{\mathop{\hbox{\rm cth}}}
\def\sh{\mathop{\hbox{\rm sh}}}
\def\max{\mathop{\hbox{\rm max}}}
\def\ind{\mathop{\hbox{\rm ind}}}
\def\Tr{\mathop{\hbox{\rm Tr}}}
\def\div{\mathop{\hbox{\rm div}}}
\def\Div{\mathop{\hbox{\rm Div}}}
\def\sn{\hbox{\rm sn}}
\def\cn{\hbox{\rm cn}}
\def\dn{\hbox{\rm dn}}
\def\sech{\hbox{\rm sech}}

\def\t{{\rm t}}

\def\idtens{\mathop{\hbox{\bf I}}}

\def\VelVec{\vec{v}}
\def\p{\partial}
\def\pdedom{\mathbb{R}^3}
\def\dV{dV}
\def\dS{d\vec{S}}
\def\dl{d\vec{\ell}}

\def\V{{\mathcal V}}
\def\S{{\mathcal S}}
\def\C{{\mathcal C}}
\def\A{{\mathcal A}}

\def\Div{\mathop{\hbox{\rm Div}}}
\def\Curl{\mathop{\hbox{\rm Curl}}}
\def\Grad{\mathop{\hbox{\rm Grad}}}

\newcommand{\vspacebefore}{\raisebox{0ex}[2.5ex][0ex]{\null}}

\newtheorem{theorem}{Theorem}
\newtheorem{lemma}{Lemma}
\newtheorem{proposition}{Proposition}
\newtheorem{corollary}{Corollary}

{\theoremstyle{definition} \newtheorem{definition}{Definition}
\newtheorem{example}{Example}
\newtheorem{remark}{Remark}
}

\renewcommand{\thefootnote}{\alph{footnote})}

\title{Local and global conservation laws of partial differential equations\\ on moving domains in three dimensions}

\author{ \renewcommand{\thefootnote}{\alph{footnote}}
  Alexey F. Shevyakov\footnotemark[1]\vspace{0.5cm},~~
  Stephen C. Anco\footnotemark[2]\\
\small $^{\rm 1}$\emph{Department of Mathematics and Statistics, University of Saskatchewan, Saskatoon, Canada}\vspace{0.2cm}\\
\small $^{\rm 2}$\emph{Department of Mathematics and Statistics, Brock University, St. Catharines, Canada}\vspace{0.2cm}\\
 }

\footnotetext[1]{Corresponding author. Previously published as Alexei F.\ Cheviakov. Electronic mail: shevyakov@math.usask.ca}
\footnotetext[2]{Electronic mail: sanco@brocku.ca}

\date{October 30, 2026}
\maketitle
\numberwithin{equation}{section}
\numberwithin{remark}{section}
\numberwithin{lemma}{section}
\numberwithin{proposition}{section}

\begin{abstract}
For PDE systems in three spatial dimensions,
firstly, a review of the various types of global conservation laws
--- volumetric, surface-flux, circulatory ---
is given using vector calculus, 
and their corresponding local forms
--- divergence, curl, gradient --- 
are explained. 
Secondly,
an extension of each type of conservation law to moving domains
is presented,
where the domain --- volume, surface, curve --- is transported by a prescribed velocity field.
This provides a generalization of Reynolds' transport theorem known for material domains
in which the velocity field is the fluid velocity. 
The main results are illustrated by numerous physical examples taken from
gas dynamics, fluid flow, nonlinear acoustics, chemotaxis, plasma dynamics,
and also Maxwell's equations, which do not involve any physically specified velocity field.
As an additional main result,
an extension of the multiplier method for finding local conservation laws of divergence-type
is presented for local curl-type and gradient-type conservation laws
on both static and moving domains. 
\end{abstract}

\section{Introduction}\label{sec:intro}

A local conservation law of a system of partial differential equations (PDEs)
is a divergence expression that vanishes for all solutions of the system.
In the case of evolution PDEs and other dynamical PDE systems,
the divergence takes the form of a continuity equation 
given by the time derivative of a density plus the spatial divergence of a flux,
formulated on a static spatial domain.
In the case of spatial PDE systems, such as elliptic PDEs, 
a local conservation law on a fixed spatial domain is purely spatial
and represents a conserved flux.
This latter type of conservation law can also exist for dynamical PDE systems
when spatial constraints are present. 

For a given PDE system,
local conservation laws have many important uses. 
They serve as mathematical expressions for fundamental physical principles,
such as conservation of mass, momentum, charge, and energy
(e.g. \cite{olver1993applications, bluman2010applications}).
They are widely employed in analysis, in particular,
studies of existence, uniqueness and stability of solutions of nonlinear PDEs
(e.g. \cite{Str-book,Gey.Pel-book}).
  
They provide a means for introducing potential variables and other nonlocal variables,
including physical scalar and vector potentials and stream functions \cite{bluman2010applications}. 
In the case of linear PDEs, they enable seeking exact solutions through, for example, the Fokas method \cite{deconinck2014method, colbrook2019unified}.
They also allow recasting a PDE system into a divergence form required by
numerical solvers, such as finite difference, finite element, and finite volume methods \cite{leveque1992numerical}.
From a mathematical perspective, local conservation laws are
an intrinsic coordinate-free aspect of the structure of a PDE system. 

In three spatial dimensions,
for dynamical PDE systems, 
the integral of a local conservation law
over any fixed spatial volume (with no source or sink terms)
gives a balance equation stating that
the time rate of change of some total quantity in the volume 
is equal to the flux through the boundary of the spatial volume.
For spatial PDE systems,
the corresponding balance equation is simply a vanishing flux integral
that represents a topological quantity. 
This type of topological conservation law can also arise from spatial constraints
in a dynamical PDE system. 

In addition to these standard divergence-type conservation laws,
there are two other important types of conservation laws in three dimensions
\cite{Anc.Che.2020}:
curl-type and gradient-type.

A local conservation law of curl type is given by 
an antisymmetric matrix expression whose divergence vanishes
for all solutions of the PDE system.
For dynamical PDE systems,
this divergence consists of the time derivative of a vector-flux density 
plus the spatial curl of a circulation vector. 
Its integral over any fixed spatial surface 
gives a balance equation whose form depends on whether the surface is closed or open.
In the case of an open surface,
the balance equation states the time rate of change of the total flux through the surface
is equal to a circulation current on the boundary of the surface.
This boundary current is absent in the case of a closed surface,
so that the balance equation gives a time-independent flux through the surface.
For a spatial PDE system,
only the circulation vector appears in the matrix-divergence expression,
and the balance equation reduces to the circulation current if the surface is open, 
or vanishes if the surface is closed. 

A local conservation law of gradient-type consists of
a totally-antisymmetric tensor expression whose divergence vanishes
for all solutions of the PDE system.
This tensor divergence in the case of dynamical PDE systems
is given by the time derivative of a vector density 
plus the spatial gradient of a flux scalar. 
The corresponding balance equation is formulated on any fixed spatial curve,
and states that the time derivative of the total circulation is equal to a scalar flux
(in the tangent direction) at the endpoints of the curve.
In the case of a closed curve,
the endpoint fluxes are absent whereby the balance equation yields
a time-independent circulation along the curve.
For a spatial PDE system,
the tensor-divergence expression contains only the scalar flux,
and thus the balance equation reduces to the endpoint fluxes,
or simply vanishes if the curve is closed.

The primary purpose of the present work is to discuss a generalization of 
all three types of local and global conservation laws
to the situation when the integration domain --- volume, surface, or curve ---
is moving with a flow given by a specified velocity field.
For conservation laws formulated on a moving volume, 
this generalization is well known in fluid mechanics
as Reynolds' transport theorem \cite{Rey}.
It was derived originally for material volumes transported by fluid flow
where the velocity field is the fluid velocity. 
Counterparts for moving surfaces and curves are much less familiar,
but are known in their global form in the literature on fluid mechanics 
\cite{McDon.Wit.1984,Gur.Str.Will.1989,Gra.1995,Cer.Fri.Gur.2005,Anc.2013,Niv.Cor.Kai.Sch.Nao}
and include Kelvin's circulation theorem as one example.

As main results,
first, 
a simple unified derivation of the three different types of moving conservation laws
will be given using vector calculus in three dimensions, 
with the velocity field being arbitrary.
Both the local and global forms of these conservation laws will be treated,
in contrast to the focus on global formulations seen in the fluid mechanics literature.
The local formulation will utilize the coordinate-free $\div$, $\curl$, $\grad$ operators.
Second, each type of moving conservation law will be illustrated by considering
a variety of physical PDE systems:
\begin{itemize}
\item
Maxwell's electromagnetic equations in vacuum;
\item
Euler equations of compressible adiabatic gas dynamics;
\item
Keller--Segel model of chemotaxis \cite{keller1971model};
\item
Navier-Stokes equations of incompressible constant-density fluid flow;
\item
nonlinear acoustic equations \cite{Ham.Bla.1998};
\item
equations of incompressible magnetohydrodynamics (MHD).
\end{itemize}

It is well known that all local divergence-type conservation laws admitted by
a given PDE system can be obtained systematically by a general method based on multipliers  
(see \cite{olver1993applications, KrasVinBook1999, bluman2010applications}
and references therein).
The multiplier method generalizes Noether's first theorem \cite{martinez1979noether},
being applicable to any PDE system whether or not it is variational. 
This method will be shown to be equally applicable
for finding divergence-type conservation laws on moving volumes. 
As a further main result,
an extension of the multiplier method will be presented
which systematically yields all of the local curl-type and gradient-type conservation laws,
on both static and moving surfaces and curves, 
admitted by a given PDE system,
under the regularity assumptions stated in Section~\ref{sec:multipliers}. 

In Section~\ref{sec:PDEsystems},
the PDE systems listed above will be reviewed using a unified notation. 
Divergence-type, curl-type, and gradient-type conservation laws on static domains
will be reviewed in Section~\ref{sec:CL:static}. 
Spatial conservation laws of topological type,
which involve no time derivatives and arise from spatial constraints,
are discussed in Section~\ref{sec:CL:topol}. 
In Section~\ref{sec:CL:moving},
their generalization to moving domains will be formulated
in both global and local forms
and illustrated by numerous physical examples.
In Section \ref{sec:multipliers},
the origin of all three types of conservation laws in local form
from multipliers is presented,
which allows for a completely systematic computation of conservation laws
for dynamical PDE systems. 
Some concluding remarks are made
in Section \ref{sec:conclusion}. 

The underlying mathematical concept for conservation laws on moving domains
comes from the Lie derivative,
which allows for a wider generalization to
moving $p$-dimensional submanifolds in $n$-dimensional space \cite{Anc.2013,Niv.Cor.Kai.Sch.Nao}
by means of the Leibniz rule. 
A vector-calculus proof of the latter in three dimensions is presented in Appendix~\ref{sec:proof},
and a Lie-derivative proof is provided in Appendix~\ref{sec:liederivative}. 
Expressions for the $\div$, $\curl$, $\grad$ operators
in an arbitrary (curvilinear) coordinate system
are summarized in Appendix~\ref{sec:curvilin.coords}.

\section{Running examples of PDE systems}\label{sec:PDEsystems}

The PDE systems listed in Section~\ref{sec:intro} will be summarized here
in a unified vector-calculus notation.
Constants are specific to each model
(for instance, $\mu$ denotes the magnetic permeability in the MHD model
and a viscosity coefficient in the acoustics model).

\subsection{Maxwell's electromagnetic equations in vacuum}\label{sec:Max:vac}

\begin{subequations}\label{eq:Max:vac}
\begin{align}
    &  \vec{E}_t  - c\, \curl \vec{B} =0, \label{eq:Max:Edyn}\\[1ex]
    &  \vec{B}_t  + c\, \curl \vec{E} =0, \label{eq:Max:Bdyn}\\[1ex]
    &  \div \vec{E} = 0, \label{eq:Max:divE} \\[1ex]
    &  \div \vec{B} = 0,\label{eq:Max:divB}
\end{align}
\end{subequations}
where $\vec{E}$ denotes the electric field, $\vec{B}$ is the magnetic field,
and $c$ is the speed of light in vacuum.
We use rationalized electromagnetic units, with energy density
$\tfrac{1}{2}(|\vec{E}|^2+|\vec{B}|^2)$ and Poynting vector $c\vec{E}\times\vec{B}$.

\subsection{Euler equations of compressible gas dynamics}\label{sec:Euler.gas.entropy}

\begin{subequations}\label{eq:gas.dyn}
\begin{align}
    &  \rho_t + \div (\rho \vec{u}) =0, \label{eq:gas.dyn:conti}\\[1ex]
    &  \rho(\vec{u}_t +(\vec{u}\cdot\nabla)\vec{u})=-\grad {p}, \label{eq:gas.dyn:mom}\\[1ex]
      & p_t +(\vec{u}\cdot\nabla)p + F(\rho,p) \div \vec{u} =0, \label{eq:gas.dyn:pres}
\end{align}
\end{subequations}
where $F= -\rho\dfrac{\partial S}{\partial \rho}\Big/\dfrac{\partial S}{\partial p}$ is a general equation-of-state function, determined from the entropy $S=S(\rho,p)$.
A commonly considered polytropic equation of state is given by 
\beq\label{eq:gas.dyn:polytropic}
F = \gamma p, 
\eeq
where the entropy is $S=c_v \ln\big( (p/p_0) (\rho/\rho_0)^{-\gamma}\big)$, 
$\gamma$ is the adiabatic exponent,
$c_v$ is the specific heat at constant volume, 
and $p_0$, $\rho_0$ are constant characteristic values of $p$ and $\rho$.

\subsection{Keller--Segel model of chemotaxis}\label{sec:KS.model}

\begin{equation}\label{eq:KSiegel}
\begin{aligned}
    \rho_t & =\nu_\rho \Delta \rho - k\,\div (\rho \grad c) ,\\[1ex]
    c_t & =\nu_c \Delta c + g(\rho) - h(c),
\end{aligned}
\end{equation}
where $\rho$ and $c$ are respectively
the bacterial density and the attractant concentration,
$\nu_\rho$ and $\nu_c$ are diffusion coefficients,
$k$ denotes the chemotactic sensitivity,
$g(\rho)$ describes the rate of attractant production,
and $h(c)$ is the rate of attractant decay.

\subsection{Navier-Stokes equations of incompressible constant-density fluid flow}\label{sec:NS3d.constdens}

\begin{subequations}\label{eq:NS3d:constd:consf}
\begin{align}
    &  \div \vec{u}=0, \label{eq:NS3d:constd:consf:div}\\[1ex]
    &  \vec{u}_t +(\vec{u}\cdot\nabla)\vec{u}=-\grad {p}+ \nu \,\Delta \vec{u},\label{eq:NS3d:constd:consf:mom}
\end{align}
\end{subequations}
where $\vec{u}$ is the fluid velocity, $p$ denotes the hydrostatic pressure divided by the constant fluid density,
$\Delta$ is the Laplace operator, and $\nu$ is the kinematic viscosity coefficient.

An alternative way to write the velocity equation \eqref{eq:NS3d:constd:consf:mom} uses the identity
\beq\label{eq:importantID:curlcurl}
(\vec{u}\cdot\nabla)\vec{u} =
\tfrac{1}{2}\grad |\vec{u}|^2 +  (\curl \vec{u})\times \vec{u},
\eeq
which yields 
\beq\label{eq:mom:vxcuv}
\vec{u}_t + (\curl \vec{u})\times \vec{u} =-\grad {P}+ \nu \,\Delta \vec{u},
\eeq
where
\beq\label{eq:modified.pressure}
P = p +\tfrac{1}{2}|\vec{u}|^2
\eeq
is called the modified pressure.

Another reformulation involves the vorticity, 
defined as 
\beq\label{eq:vorticity}
\vec{\omega}=\curl \vec{u}
\eeq
which satisfies 
\begin{subequations}\label{eq:NS3d:constd:consf:vort}
\begin{align}
 &\div  \vec{\omega}=0, \label{eq:NS3d:constd:consf:vort:divom}\\[1ex]
 &\vec{\omega}_t+ \curl(\vec{\omega}\times\vec{u})= \nu \,\Delta \vec{\omega}\,. \label{eq:NS3d:constd:consf:vort:momom}
\end{align}
\end{subequations}
where equation \eqref{eq:NS3d:constd:consf:vort:momom} is obtained
by taking the curl of the velocity equation \eqref{eq:mom:vxcuv}.

\subsection{Nonlinear acoustic equations}\label{sec:acoustic}

\begin{subequations}\label{eq:acoustic}
  \begin{align}
    & \rho_t +\rho_0 \div \vec{u} =0, \label{eq:acoustic.rho}\\[1ex]
    & \rho_0 \vec{u}_t = - \grad p + \mu\, \Delta \vec{u} + \nu\, \grad\div \vec{u}, \label{eq:acoustic.vel} \\[1ex]
    & p = c_0^2\,\rho + f(\rho) \label{eq:acoustic.eos}
  \end{align}
\end{subequations}
where $\vec{u}$ is the air velocity, $\rho$ is the density change in the disturbed air,
$p$ is the pressure change in the disturbed air,
$\rho_0$ is the density of the undisturbed air,
and $c_0$ is speed of sound in the undisturbed air,
while $\mu$ and $\nu$ are viscosity coefficients. 
The equation of state \eqref{eq:acoustic.eos} involves a nonlinear function $f(\rho)$ with $f(0)=f'(0)=0$,
so that $c_0$ is the sound speed of the linearized model. 
This function is typically approximated by the lowest-order term in a Taylor expansion:
$f(\rho) = \beta \rho^2$, which yields
\begin{equation}\label{eq:acoustic.p}
  p = c_0^2\rho + \beta \rho^2 .
\end{equation}

Taking the divergence of the second PDE \eqref{eq:acoustic.vel} and substituting $\div\vec{u}$ from the first PDE \eqref{eq:acoustic.rho}
leads to a nonlinear wave equation for the density:
\beq\label{eq:acoustic.dens.waveeqn}
\rho_{tt} = c_0^2 \Delta \rho + \Delta\left( f(\rho) +\dfrac{\mu+\nu}{\rho_0} \rho_t \right).
\eeq
In the case when the approximation \eqref{eq:acoustic.p} is considered,
inverting it gives the approximation
$\rho = c_0^{-2} \left( p - \beta c_0^{-4} p^2 \right)$
holding to quadratic order in $p$ for small $p$, 
whereby the nonlinear wave equation becomes
\beq\label{eq:acoustic.pres.waveeqn}
\left( p-c_0^{-4}\beta p^2 \right)_{tt} \simeq c_0^2 \Delta p + \dfrac{\mu+\nu}{\rho_0} \Delta \left( p - \beta c_0^{-4} p^2 \right)_t .
\eeq

The curl of the PDE \eqref{eq:acoustic.vel} yields
\beq\label{eq:acoustic.vort}
 \vec{\omega}_t = \frac{\mu}{\rho_0} \,\Delta \vec{\omega}
\eeq
for the vorticity \eqref{eq:vorticity}.
This is a vector heat equation.

\subsection{Equations of incompressible magnetohydrodynamics}\label{sec:MHD.full}
    
\begin{subequations}\label{eq:MHD:full}
\begin{align}
    &  \rho_t + \div (\rho\vec{u}) = 0,\label{eq:MHD:full:conti}\\[1ex]
    &
       \vec{u}_t = \vec{u}\times \curl \vec{u} + \dfrac{1}{\mu\rho}\curl \vec{B} \times \vec{B} -\dfrac{1}{\rho}\grad P 
       - \frac{1}{2} \grad |\vec{u}|^2 + \dfrac{\lambda}{\rho} \Delta \vec{u},
       \label{eq:MHD:full:mom}\\[1ex]
    &  \vec{B}_t = \curl(\vec{u} \times \vec{B}) + \eta\, \Delta \vec{B},\label{eq:MHD:full:Bdynamics}\\[1ex]
    &  \div \vec{B} = 0, \label{eq:MHD:full:divB}\\[1ex]
    &  \div \vec{u} = 0,\label{eq:MHD:full:divV}
\end{align}
\end{subequations}
where $\rho>0$ denotes plasma density, $\vec{B}$ is the magnetic field,
$\vec{J} = \dfrac{1}{\mu}\curl \vec{B}$ is the electric current density,
$\vec{u}$ is the plasma velocity, $P$ denotes pressure,
$\mu$ is the constant magnetic permeability coefficient,
$\eta$ is the constant magnetic diffusivity,
and $\lambda$ is a constant dynamic viscosity coefficient.

\section{Dynamical conservation laws on static domains}\label{sec:CL:static}

For a general dynamical PDE system in three dimensions,
the three types of conservation laws \cite{Anc.Che.2020}
will now be reviewed in their local and global forms.
Let $\depvars$ denote the dependent variable(s) in the system,
which are functions of time $t$ and position $\vec{x}$.
The notation $f[\depvars]$ will be used to denote a differential function 
that depends on $\depvars$ and partial derivatives of $\depvars$
with respect to $t$ and $\vec{x}$ up to some finite order,
as well as explicitly on $t$ and $\vec{x}$ themselves. 
Likewise, $\vec{f}[\depvars]$ will denote a differential vector function,
and $\underline{f}[\depvars]$ will denote a differential matrix function. 
In the examples below,
Cartesian coordinates are denoted $\vec{x}=(x^1,x^2,x^3)=(x,y,z)$,
with the corresponding unit basis vectors
$\hat{\vec{e}}_1=\hat{\vec{x}}$, $\hat{\vec{e}}_2=\hat{\vec{y}}$, $\hat{\vec{e}}_3=\hat{\vec{z}}$;
$\epsilon^{ijk}$ denotes the Levi-Civita symbol,
and summation over repeated indices is assumed.

\subsection{Volumetric conservation laws} \label{sec:vol:staticdomain}

\begin{figure}[htbp]
\begin{subfigure}{0.3\textwidth}
 \centering\includegraphics[width=0.9\textwidth]{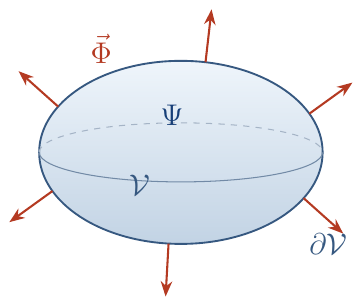}
 \caption{}\label{fig:schematics_3D}
\end{subfigure}
\hfill
\begin{subfigure}{0.3\textwidth}
\centering\includegraphics[width=0.95\textwidth]{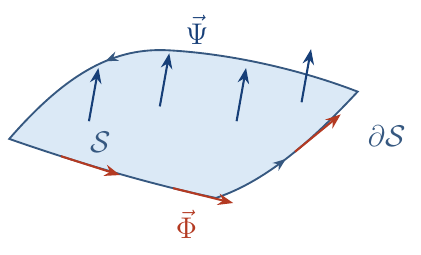}
\caption{}\label{fig:schematics_2D}
\end{subfigure}
\hfill
\begin{subfigure}{0.3\textwidth}
\centering\includegraphics[width=0.95\textwidth]{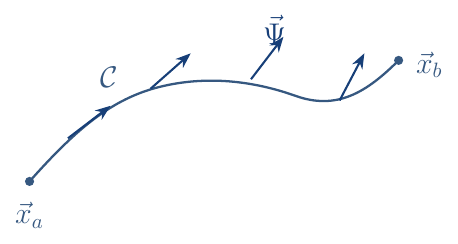}
\caption{}\label{fig:schematics_1D}
\end{subfigure}
\hfill
\caption{Global conservation laws: (\subref{fig:schematics_3D}) volumetric, (\subref{fig:schematics_2D}) surface-flux, (\subref{fig:schematics_1D}) circulatory.}
\label{fig:schematics}
\end{figure}

A \emph{global volumetric conservation law} of a dynamical PDE system 
is a balance equation 
\begin{equation}\label{3Dglobalvolconslaw}
\frac{d}{dt} \int_{\V} \ConsLawDensity\, \dV = -\oint_{\p\V} \vec{\ConsLawFlux} \cdot \,\dS\,
\end{equation}
holding for all solutions $\depvars(t, \vec{x})$ of the system. 
Here $\V \subseteq\pdedom$ is a closed volume having boundary surface $\p\V$,
which is assumed to be piecewise smooth, 
with $\dS=\hat{\vec{n}}\,dS$ being the surface element
given in terms of the outward unit normal  vector $\hat{\vec{n}}$
and the area element $dS$. 
See Figure \ref{fig:schematics_3D}. 
The scalar differential function $\ConsLawDensity=\ConsLawDensity[\depvars]$ is called the \emph{conserved density},
and the vector differential function $\vec{\ConsLawFlux}=\vec{\ConsLawFlux}[\depvars]$ is called
the \emph{spatial flux} for the conservation law \eqref{3Dglobalvolconslaw}.

The physical meaning of a global conservation law \eqref{3Dglobalvolconslaw}
is that the rate of change of the volumetric quantity
\begin{equation}\label{3DvolC}
C[\depvars;\V]= \int_{\V}  \ConsLawDensity\,\dV
\end{equation}
in any volume $\V\subseteq\pdedom$ is balanced by the net surface flux
\begin{equation}\label{3Dvolflux}
F[\depvars;\p\V] = \oint_{\p\V} \vec{\ConsLawFlux}\cdot \,\dS
\end{equation}
escaping through the volume boundary surface $\p\V$.
In particular,
this balance can be interpreted as an absence of sources or sinks in the volume $\V$.

A volumetric conservation law \eqref{3Dglobalvolconslaw} that holds for all volumes $\V$ in $\pdedom$
is equivalent, through the divergence theorem,
to a local continuity equation
\begin{equation}\label{3Dvolconslaw}
\partial_ t\ConsLawDensity + \div \vec{\ConsLawFlux} =0
\end{equation}
satisfied for all solutions $\depvars(t, \vec{x})$ of the given PDE system.

\medskip\noindent\textbf{Example 1.}
For the Navier-Stokes equations of incompressible, constant-density fluid flow \eqref{eq:NS3d:constd:consf},
the following volumetric conservation laws hold.

\begin{enumerate}
\item
Generalized mechanical momentum,
which comprises three independent conservation laws
in the $x$, $y$, $z$ directions
\beq\label{eq:NS:CL:genMom}
\ConsLawDensity^{(i)} = f(t)\, u^i,
\qquad
\vec{\ConsLawFlux}^{(i)} = f(t) \left( u^i\, \vec{u} - \nu\,\grad u^i + p\,\hat{\vec{e}}_i \right) - f'(t)\, x^i\, \vec{u},
\qquad i=1,2,3,
\eeq
involving an arbitrary function $f(t)$.
These expressions can be written in a vector-matrix form as 
\beq\label{eq:NS:CL:genMom:vec}
\vec{\ConsLawDensity} = f(t)\vec{u},
\qquad
\underline{\ConsLawFlux} = f(t)\big( \vec{u}\vec{u}   -\nu (\grad \vec{u})^\t + p\idtens   \big) - f'(t) \vec{u}\, \vec{x}
\eeq
where $\t$ denotes transpose and $\idtens$ is the identity matrix.\footnote{Note  $[\div(\vec{u}\vec{x})]^i=\sum_j\partial_{x^j}(u^j x^i)=x^i\div\vec{u}+u^i=u^i$, since $\div\vec{u}=0$.}
Our Cartesian component conventions are
$(\vec{a}\vec{b})^{ij}=a^i b^j$ and
$(\grad\vec{a})^{ij}=\partial_{x^j}a^i$.
Then the local conservation laws \eqref{eq:NS:CL:genMom} take the form
\begin{equation}
  \partial_t \vec{\ConsLawDensity} + \div \underline{\ConsLawFlux} =0
\end{equation}
where 
$\div\underline{A}$ denotes the vector with $x-$component
$\partial_ x \underline{A}^{11} + \partial_ y \underline{A}^{21} + \partial_ z \underline{A}^{31}$,
and similarly for $y-$ and $z-$ components,
for any differential matrix function $\underline{A}$.

\item
  Angular momentum,
  which comprises three independent conservation laws
  in the $x$, $y$, $z$ directions
\beq
\ConsLawDensity^{(i)} = \epsilon^{ijk} x^j u^k,
\qquad
\vec{\ConsLawFlux}^{(i)} = \epsilon^{ijk} x^j u^k\, \vec{u}
- \nu\,\grad \left( \epsilon^{ijk} x^j u^k \right)
+ \epsilon^{ijk} \left( p\, x^j - 2\nu\, u^j \right) \hat{\vec{e}}_k
\eeq
These can be written succinctly as 
\beq\label{eq:NS:CL:angMom}
\vec{\ConsLawDensity} = \vec{x}\times\vec{u},
\qquad
\underline{\ConsLawFlux} = \vec{u} (\vec{x}\times\vec{u})
-\nu(\grad (\vec{x}\times\vec{u}))^\t
-2\nu \underline{\vec{u}} + p\underline{\vec{x}}
\eeq
where 
$\underline{\vec{A}}$ denotes an antisymmetric matrix with components
$\underline{\mathrm{A}}^{12} = -\underline{\mathrm{A}}^{21} = A^3$,
$\underline{\mathrm{A}}^{23} = -\underline{\mathrm{A}}^{32} = A^1$,
$\underline{\mathrm{A}}^{31} =-\underline{\mathrm{A}}^{13} =A^2$,
for any differential vector function $\vec{A}$
(for details, see \cite{cheviakov2014generalized}).
This matrix satisfies the identity
\beq\label{eq:div.underline}
\div\, \underline{\vec{A}} = -\curl \vec{A} ,
\eeq
which will be used below.
Together with the generalized continuity equation
discussed in Section~\ref{sec:CL:topol},
above set of divergence-type conservation laws
of the Navier-Stokes equations ($\nu \ne 0$) is complete \cite{gusyatnikova1989symmetries}.

\item
 In the case of zero viscosity, $\nu=0$,
kinetic energy 
\beq\label{eq:Eul:CL:KinEn}
\ConsLawDensity = K = \tfrac{1}{2} |\vec{u}|^2 ,
\quad
\ConsLawFlux = (K+p) \vec{u} , 
\eeq
and helicity
\beq\label{eq:Eul:CL:heli}
\ConsLawDensity = h = \vec{u}\cdot \vec{\omega},
\quad
\ConsLawFlux = \vec{u}\times \nabla (\tfrac{1}{2}|\vec{u}|^2 +p) 
+ \vec{u}\times (\vec{u}\times \vec{\omega})
\eeq
which describes the knottedness of vorticity lines.

\end{enumerate}

\medskip\noindent\textbf{Example 2.}
The Euler equations of compressible gas dynamics \eqref{eq:gas.dyn}
possess the following volumetric conservation laws.
\begin{enumerate}
\item
Mass   
\beq\label{eq:gas:CL:mass}
\ConsLawDensity = \rho,
\qquad
\ConsLawFlux = \rho \vec{u}
\eeq

\item
 Entropy
\beq\label{eq:gas:CL:entropy}
\ConsLawDensity = \rho S,
\qquad
\vec{\ConsLawFlux} = \rho S \vec{u}
\eeq
More generally, $\ConsLawDensity=\rho f(S)$ and
$\vec{\ConsLawFlux}=\rho f(S)\vec{u}$ give a conservation law
for any function $f$.

\item
  Energy
\beq\label{eq:gas:CL:energy}
\ConsLawDensity = E= \tfrac{1}{2} \rho |\vec{u}|^2 +\rho\, e,
\qquad
\ConsLawFlux = (E+p)\vec{u} 
\eeq
where
\beq\label{eq:gas:CL:internalenergy}
e = \int (P(\rho,S)/\rho^2)\, d\rho
\eeq
is the internal energy,
with $P$ being the pressure $p$ expressed as a function of $\rho$ and $S$
through inverting $S=S(\rho,p)$. 
For the polytropic equation of state with $\gamma\ne1$,
$e = \frac{1}{\gamma-1} p/\rho$.

\item
  Momentum, comprising three conservation laws ($i=1,2,3$),
\beq
\ConsLawDensity^{(i)} = \rho\, u^i,
\qquad
\vec{\ConsLawFlux}^{(i)} = \rho\, u^i\, \vec{u} + p\,\hat{\vec{e}}_i
\eeq

\item
  Galilean momentum, comprising three conservation laws,
\beq
\ConsLawDensity^{(i)} = \rho \left( t\, u^i - x^i \right),
\qquad
\vec{\ConsLawFlux}^{(i)} = \rho \left( t\, u^i - x^i \right) \vec{u} + t\, p\,\hat{\vec{e}}_i
\eeq

\item
  Angular momentum, comprising three conservation laws,
\beq
\ConsLawDensity^{(i)} = \rho\, \epsilon^{ijk} x^j u^k,
\qquad
\vec{\ConsLawFlux}^{(i)} = \rho\, \epsilon^{ijk} x^j u^k\, \vec{u} + p\, \epsilon^{ijk} x^j\, \hat{\vec{e}}_k
\eeq

\end{enumerate}

\medskip\noindent\textbf{Example 3.}
The incompressible MHD equations \eqref{eq:MHD:full} possess the following volumetric conservation laws. 
\begin{enumerate}
  \item
Each of the three components of the right hand side of the magnetic field equation \eqref{eq:MHD:full:Bdynamics}
is a divergence expression,
which yields three conservation laws
\beq
\ConsLawDensity^{(i)} = B^i,
\qquad
\vec{\ConsLawFlux}^{(i)} = \left( u^k B^i - u^i B^k \right) \hat{\vec{e}}_k - \eta\, \grad B^i ,
\qquad i=1,2,3.
\eeq
These have the succinct vector-matrix form
\beq\label{eq:MHD:CL:magnetic}
\vec{\ConsLawDensity} = \vec{B},
\qquad
\underline{\ConsLawFlux} = \underline{\vec{u} \times \vec{B}} - \eta (\grad \vec{B})^\t 
\eeq
where $\underline{\vec{u} \times \vec{B}}$ is the antisymmetric matrix
associated with the vector $\vec{u}\times\vec{B}$, explicitly
\[
\underline{\vec{u}\times\vec{B}}=\vec{u}\vec{B}-\vec{B}\vec{u}
=\begin{pmatrix}
0 & u^1 B^2-u^2 B^1 & u^1 B^3-u^3 B^1\\
u^2 B^1-u^1 B^2 & 0 & u^2 B^3-u^3 B^2\\
u^3 B^1-u^1 B^3 & u^3 B^2-u^2 B^3 & 0
\end{pmatrix}.
\]
By the identity \eqref{eq:div.underline},
this conservation law is equivalent to
the magnetic field equation \eqref{eq:MHD:full:Bdynamics}.

\item
Momentum,
in the succinct vector-matrix notation of Section~\ref{sec:vol:staticdomain},
\beq\label{eq:MHD:CL:mom}
\vec{\ConsLawDensity} = \rho\vec{u},
\qquad
\underline{\ConsLawFlux} = \rho\, \vec{u}\vec{u}
+ \Big( P + \frac{1}{2\mu}|\vec{B}|^2 \Big) \idtens
- \frac{1}{\mu} \vec{B}\vec{B}
- \lambda\, (\grad\vec{u})^\t
\eeq
where the magnetic terms are the negative of the magnetic Maxwell stress tensor in the plasma,
obtained through the identity
$(\curl\vec{B})\times\vec{B} = \div(\vec{B}\vec{B}) - \tfrac{1}{2}\grad|\vec{B}|^2$
holding by virtue of \eqref{eq:MHD:full:divB}.

\item
Angular momentum, with components
\beq\label{eq:MHD:CL:angmom}
\vec{\ConsLawDensity} = \vec{x}\times(\rho\vec{u}),
\qquad
\underline{\ConsLawFlux}^{\,ji} = \varepsilon^{ikl} x^k \big( \rho\, u^j u^l + \sigma^{jl} \big)
- \lambda\, \partial_{x^j} (\vec{x}\times\vec{u})^i - 2\lambda\, \underline{\vec{u}}^{\,ji}
\eeq
where
$\sigma^{jl} = \big(P + \tfrac{1}{2\mu}|\vec{B}|^2\big)\delta^{jl} - \tfrac{1}{\mu}B^jB^l$
denotes the symmetric (pressure and magnetic) part of the momentum flux,
and $\underline{\vec{u}}$ is the antisymmetric matrix associated to $\vec{u}$
as in Section~\ref{sec:vol:staticdomain}.

\end{enumerate}

\medskip\noindent\textbf{Example 4.}
In the Keller--Segel chemotaxis model \eqref{eq:KSiegel},
the first PDE constitutes a volumetric conservation law
\beq\label{eq:KSiegel:VolCL}
\ConsLawDensity = \rho,
\qquad
\ConsLawFlux = {-}\nu_\rho \grad\rho +k \rho \grad c
\eeq
This describes local conservation of the bacterial mass.

\medskip\noindent\textbf{Example 5.}
Maxwell's equations in vacuum \eqref{eq:Max:vac}
possess the following volumetric conservation laws,
arising from the spacetime symmetries of the system
(time and space translations, rotations, and Lorentz boosts).
\begin{enumerate}
\item
Energy
\beq\label{eq:Max:CL:energy}
\ConsLawDensity = \tfrac{1}{2}\big(|\vec{E}|^2+|\vec{B}|^2\big),
\qquad
\vec{\ConsLawFlux} = c\, \vec{E}\times\vec{B}
\eeq
where the flux is the Poynting vector.

\item
Momentum
\beq\label{eq:Max:CL:mom}
\vec{\ConsLawDensity} = \vec{E}\times\vec{B},
\qquad
\underline{\ConsLawFlux} = c\, \underline{T},
\qquad
\underline{T} = \tfrac{1}{2}\big(|\vec{E}|^2+|\vec{B}|^2\big)\idtens
- \vec{E}\vec{E} - \vec{B}\vec{B} ,
\eeq
where $\underline{T}$ is the negative of the Maxwell stress tensor.

\item
Angular momentum, with components
\beq\label{eq:Max:CL:angmom}
\vec{\ConsLawDensity} = \vec{x}\times(\vec{E}\times\vec{B}),
\qquad
\underline{\ConsLawFlux}^{\,ji} = c\,\varepsilon^{ikl}\, x^k\, T^{jl}
\eeq
using the symmetry of $\underline{T}$.

\item
Boost momentum
\beq\label{eq:Max:CL:boost}
\vec{\ConsLawDensity} = c\,t\, \vec{E}\times\vec{B}
- \tfrac{1}{2}\vec{x}\,\big(|\vec{E}|^2+|\vec{B}|^2\big),
\qquad
\underline{\ConsLawFlux}^{\,ji} = c^2 t\, T^{ji} - c\, x^i\, (\vec{E}\times\vec{B})^j
\eeq

\end{enumerate}

\medskip\noindent\textbf{Example  6.}
In the nonlinear acoustic equations \eqref{eq:acoustic},
the first PDE is a volumetric conservation law
\beq\label{eq:acoustic:VolCL}
\ConsLawDensity = \rho,
\qquad
\ConsLawFlux = \rho_0 \vec{u}, 
\eeq
for the excess mass of the acoustic disturbance.
Each component of the vector PDE \eqref{eq:acoustic.vel} also constitutes
a volumetric conservation law
\beq
\ConsLawDensity^{(i)} = \rho_0\, u^i,
\qquad
\vec{\ConsLawFlux}^{(i)} = \left( p - \nu\, \div \vec{u} \right) \hat{\vec{e}}_i - \mu\, \grad u^i ,
\qquad i=1,2,3,
\eeq
which describes the momentum of the moving air in the $x$, $y$, $z$ directions.

Three additional volumetric conservation laws comprise angular momentum
in $x$, $y$, $z$ directions
\beq
\ConsLawDensity^{(i)} = \rho_0\, \epsilon^{ijk} x^j u^k,
\qquad
\vec{\ConsLawFlux}^{(i)} = -\mu\, \grad \left( \epsilon^{ijk} x^j u^k \right)
+ \epsilon^{ijk} \left( \left( p - \nu\,\div\vec{u} \right) x^j - 2\mu\, u^j \right) \hat{\vec{e}}_k
\eeq

In the case $\mu=0$, there is a volumetric conservation law describing helicity:
\beq\label{eq:acoustic:heli}
\ConsLawDensity = h = \vec{u}\cdot (\curl \vec{u}),
\quad
\ConsLawFlux = \frac{1}{\rho_0} (p -\nu\,\div\vec{u}) \curl \vec{u} 
\eeq

When both viscosities vanish, $\mu=\nu=0$,
there is also a volumetric conservation law of energy:
\beq
\ConsLawDensity = \tfrac{1}{2}\rho_0 |\vec{u}|^2 + \frac{1}{\rho_0}\Big( \tfrac{1}{2}c_0^2\rho^2 + F(\rho) \Big), 
\quad
\vec{\ConsLawFlux} =p\,\vec{u}
\eeq
where $F'(\rho)=f(\rho)$, $F(0)=0$, so that the density comprises the kinetic and compression energy of the disturbance.

In addition,
three conservation laws describe the Galilean (centre-of-mass) momentum
in the $x$, $y$, $z$ directions:
\beq
\ConsLawDensity^{(i)} = \rho_0\, t\, u^i - x^i \rho,
\qquad
\vec{\ConsLawFlux}^{(i)} = t \left( \left( p - \nu\,\div\vec{u} \right) \hat{\vec{e}}_i - \mu\,\grad u^i \right) - \rho_0\, x^i\, \vec{u}
\eeq
holding for arbitrary viscosities $\mu$, $\nu$ and an arbitrary equation of state \eqref{eq:acoustic.eos}.

\subsection{Surface-flux conservation laws} \label{sec:surf:staticdomain}

The second type of global conservation law for a general dynamical PDE system is
a surface-flux balance equation \cite{Anc.Che.2020}
\begin{equation}\label{3Dglobalsurfconslaw}
\frac{d}{dt} \int_{\S} \vec{\ConsLawDensity}\cdot\dS = -\oint_{\p\S} \vec{\ConsLawFlux}\cdot d\vec{\ell}\,,
\end{equation}
where $\S$ is a fixed orientable surface in $\pdedom$ having a piecewise-smooth boundary curve $\p\S$,
with the line element $d\vec{\ell} = \hat{\vec{\ell}}\,ds$ being given in terms of
the arclength $ds$ and the unit tangent vector $\hat{\vec{\ell}}$ along this curve,
with the surface element $\dS= \hat{\vec{n}}\,dA$ oriented such that
$\hat{\vec{\ell}}\times\hat{\vec{n}}$ is pointing outward (Figure \ref{fig:schematics_2D}).
Here $\vec{\ConsLawDensity}=\vec{\ConsLawDensity}[\depvars]$ and $\vec{\ConsLawFlux}=\vec{\ConsLawFlux}[\depvars]$
are vector differential functions,
which are referred to as the \emph{conserved flux density} and the \emph{spatial circulation flux}.

A flux equation \eqref{3Dglobalsurfconslaw} holding for all solutions $\depvars(t, \vec{x})$ of the PDE system
is called a \emph{global surface-flux conservation law}.
Its physical meaning is that the rate of change of the surface integral quantity
\begin{equation}\label{3DsurfC}
C[\depvars;\S]= \int_{\S} \vec{\ConsLawDensity} \cdot \dS
\end{equation}
is balanced by the net flux circulation
\begin{equation}\label{3Dsurfflux}
F[\depvars;\p\S] = \oint_{\p\S} \vec{\ConsLawFlux} \cdot d \vec{\ell}
\end{equation}
which is trapped around the surface boundary.
When $\S$ is a closed surface having no boundary (i.e., $\p\S=\emptyset$),
such as for example a surface topologically equivalent to a sphere or to a torus,
the quantity \eqref{3DsurfC} yields a time-independent surface flux.
Conservation in this case can be understood as a consequence of the flux having nowhere to be trapped.

The global surface-flux conservation law \eqref{3Dglobalsurfconslaw} that holds for all surfaces $\S\subset \mathbb{R}^3$ has a formulation as a local vector continuity equation.
First, using the Stokes' theorem,
the surface boundary integral \eqref{3Dsurfflux} can be converted into a surface integral:
$\oint_{\p\S} \vec{\ConsLawFlux}\cdot d\vec{\ell} = \int_{\S} (\curl\vec{\ConsLawFlux})\cdot\dS$.
Then the surface integral can be combined with the surface flux integral \eqref{3DsurfC}
to get
$\int_{\S} (\partial_t\vec{\ConsLawDensity}+ \curl \vec{\ConsLawFlux})\cdot \dS =0$.
This equation will hold for an arbitrary surface $\S\subset\pdedom$ if and only if the density and boundary flux pair satisfy the local equation \cite{Anc.Che.2020}
\begin{equation}\label{3Dsurfconslaw}
\partial_ t\vec{\ConsLawDensity} + \curl\,\vec{\ConsLawFlux} =0,
\end{equation}
holding for all solutions $\depvars(t,\vec{x})$ of the given PDE system.
This type of local vector continuity equation
is called a \emph{local surface-flux conservation law}.

For a local surface-flux conservation law \eqref{3Dsurfconslaw},
it is easy to see that $\div \partial_t \vec{\ConsLawDensity} = 0$ holds for all solutions of the given PDE system. 
Therefore, if $\vec{\ConsLawDensity}$ is divergence-free at the initial time $t=0$,
it remains divergence-free for all times.

\medskip\noindent\textbf{Example 1.}
In the case of inviscid flow $\nu=0$,
the vorticity evolution equation \eqref{eq:NS3d:constd:consf:vort:momom}
\beq
\vec{\omega}_t+ \curl(\vec{\omega}\times\vec{u})= 0
\eeq
is a local surface-flux conservation law as it stands.
It corresponds to the conservation of the vorticity flux
\beq
C= \int_{\S} \vec{\omega} \cdot \dS
\eeq
whose rate of change is balanced by the circulation
\beq
F = \oint_{\p\S} (\vec{\omega}\times\vec{u}) \cdot d\vec{\ell}
\eeq
for any piecewise-smooth surface $\S$ in $\mathbb{R}^3$.
For a closed surface without a boundary,
the total vorticity flux through $\S$ stays constant.

\medskip\noindent\textbf{Example 2.}
Similarly, the magnetic field dynamics \eqref{eq:MHD:full:Bdynamics}
in the case of infinitely conductive plasmas $\eta=0$ yields
a local surface-flux conservation law 
\beq
\vec{B}_t = \curl(\vec{u} \times \vec{B}) . 
\eeq
Here the global conserved quantity is the total magnetic flux
\beq\label{eq:sflux:Cb}
C_B= \int_{\S} \vec{B} \cdot \dS .
\eeq

\medskip\noindent\textbf{Example 3.}
For  Maxwell's equations in a vacuum,
both the electric and magnetic field dynamics equations \eqref{eq:Max:Edyn} and \eqref{eq:Max:Bdyn}
\beq
    \vec{E}_t  - c \curl \vec{B} =0, \qquad  \vec{B}_t  + c \curl \vec{E} =0
\eeq
constitute local surface-flux conservation laws that correspond to the global conservation of the total fluxes
\beq
C_E= \int_{\S} \vec{E} \cdot \dS,
\quad
C_B= \int_{\S} \vec{B} \cdot \dS
\eeq
through any closed surface $\S$.

\subsection{Circulatory conservation laws} \label{sec:circ:staticdomain}

A line integral equation
\begin{equation}\label{3Dglobalcurvconslaw}
  \dfrac{d}{dt} \int_{\C} \vec{\ConsLawDensity} \cdot d\vec{\ell}
  = -\ConsLawFlux\Big|_{\p\C}
\end{equation}
is the third type of global conservation law for a dynamical PDE system,
where $\C$ is a piecewise-smooth curve with endpoints $(\vec{x}_a,\vec{x}_b)=\p\C$
(see Figure \ref{fig:schematics_1D}).
Here $\vec{\ConsLawDensity} =\vec{\ConsLawDensity} [\depvars]$
is a vector differential function,
which is called the \emph{conserved circulation density},
and $\ConsLawFlux=\ConsLawFlux[\depvars]$ is a scalar differential function,
which is called the \emph{endpoint flux}.
A line integral equation \eqref{3Dglobalcurvconslaw} holding for all solutions $\depvars(t, \vec{x})$ of the PDE system
is a \emph{global circulatory conservation law} \cite{Anc.Che.2020}. 

The physical meaning of a global conservation law \eqref{3Dglobalcurvconslaw}
describes a balance between the rate of change of the line integral quantity
\begin{equation}\label{3DcurvC}
C[\depvars;\C]= \int_{\C} \vec{\ConsLawDensity}\cdot d\vec{\ell}
\end{equation}
and a net outward flow through the endpoints $(\vec{x}_a,\vec{x}_b)$ of the curve.
If $\C$ is a closed curve, $\p\C=\emptyset$,
then the total circulation \eqref{3DcurvC} is a constant of motion.

A global circulatory conservation law \eqref{3Dglobalcurvconslaw} that holds for all curves $\C\subset \pdedom$ can be formulated as a local continuity equation
by using the fundamental theorem of calculus for line integrals.
In particular,
this allows converting the flux at the endpoints on the right hand side of equation \eqref{3Dglobalcurvconslaw}
into a gradient line integral 
$\ConsLawFlux\big|_{\p\C} = \int_{\C} \grad \ConsLawFlux \cdot d\vec{\ell}$
which can be combined with the line integral on the left hand side to get
$\int_{\C} (\partial_t\vec{\ConsLawDensity}+ \grad \ConsLawFlux )\cdot d\vec{\ell} =0$. 
This equation will hold for an arbitrary curve $\C\subset\pdedom$ if and only if the density and boundary flux pair satisfy the local equation \cite{Anc.Che.2020}
\begin{equation}\label{3Dcurvconslaw}
\partial_t\vec{\ConsLawDensity} + \grad\,\ConsLawFlux =0,
\end{equation}
which will be referred to as a \emph{local circulatory conservation law},
holding for all solutions $\depvars(t, \vec{x})$ of the given PDE system. 

For circulatory conservation laws \eqref{3Dcurvconslaw},
it is easy to see that $\curl \partial_t\vec{\ConsLawDensity} = 0$ holds for all solutions of the given PDE system.
As a consequence, if $\vec{\ConsLawDensity}$ is irrotational at the initial moment of time, then it is irrotational at all times.

\medskip\noindent\textbf{Example.}
For the incompressible constant-density Euler equations,
given by the Navier-Stokes equations \eqref{eq:NS3d:constd:consf} with $\nu=0$,
Beltrami flows are defined by $\curl \vec{u} = \alpha \vec{u}$
for some scalar function $\alpha=\alpha(\vec{x},t)$.
When the velocity equation \eqref{eq:mom:vxcuv} with $\nu=0$
is specialized to a Beltrami flow,
the term $(\curl\vec{u})\times\vec{u}=\alpha\,\vec{u}\times\vec{u}$ vanishes,
which yields a local circulatory conservation law
\beq\label{eq:mom:vxcuv:circ}
\vec{u}_t + \grad P=0 
\eeq
involving the modified pressure \eqref{eq:modified.pressure}.
Its global form is given by 
\beq
\dfrac{d}{dt} \int_{\C} \vec{u} \cdot d\vec{\ell} = -P\big|_{\p\C}
\eeq
with the endpoint terms vanishing when $\C$ is a closed loop, $\p\C=\emptyset$. 
This is a version of Kelvin's circulation theorem for fixed curves in Beltrami flows.

\section{Spatial topological conservation laws}\label{sec:CL:topol}

For completeness, this section briefly reviews the spatial topological conservation laws
developed in Ref.~\cite{Anc.Che.2020},
so as to complete the classification and contrast them with the moving-domain laws below.
In addition to the dynamical conservation laws reviewed in Section~\ref{sec:CL:static},
a PDE system may admit purely spatial local conservation laws,
in which no time derivatives appear
(although the flux may depend on time parametrically),
where the spatial flux has vanishing divergence, curl, or gradient
for all solutions of the system.
Such conservation laws typically arise from spatial constraints,
and their global forms express topological information about solutions.

First, a \emph{topological flux conservation law} is a vanishing divergence expression
\cite{Anc.Che.2020}
\beq\label{eq:CL:topol:vol:loc}
\div\,\vec{\ConsLawFlux}=0
\eeq
holding for all solutions $\depvars(t,\vec{x})$ of a given PDE system.
By the divergence theorem, its global form on a domain $\V\subseteq\pdedom$,
\beq\label{eq:CL:topol:vol:glob}
\oint_{\p\V}\, \vec{\ConsLawFlux} \cdot \dS =0 ,
\eeq
expresses a balance of fluxes through the connected components of the boundary $\p\V$:
for example, when $\p\V$ consists of two connected components $\S_1$ and $\S_2$,
as for a spherical shell,
the flux entering through $\S_1$ leaves through $\S_2$
at every instant of time.

\medskip\noindent\textbf{Example 1.}
For the incompressible Navier-Stokes (or Euler) equations \eqref{eq:NS3d:constd:consf},
a generalized continuity equation
\beq\label{eq:NS:CL:genConti}
\div\,\big(k(t)\, \vec{u}\big)=0
\eeq
holds for an arbitrary function $k(t)$.
The divergence constraints \eqref{eq:MHD:full:divB} in MHD
and \eqref{eq:Max:divE}--\eqref{eq:Max:divB} in electromagnetism
yield the same kind of global flux-balance relations
for the fields $\vec{B}$ and $\vec{E}$.

Second, a \emph{topological circulation conservation law} is a vanishing curl expression
\beq\label{eq:CL:topol:surf:loc}
\curl\,\vec{\ConsLawFlux}=0 .
\eeq
By Stokes' theorem, its global form is the vanishing circulation
\beq\label{eq:CL:topol:surf:glob}
\oint_{\p\S}\, \vec{\ConsLawFlux} \cdot d\vec{\ell} = 0
\eeq
around the boundary of any orientable surface $\S\subset\pdedom$;
when the boundary $\p\S$ has several components,
equation \eqref{eq:CL:topol:surf:glob} becomes a balance of circulations.

\medskip\noindent\textbf{Example 2.}
For equilibrium configurations of ideal MHD plasmas
($\eta=\lambda=0$ and $\p_t\equiv0$ in \eqref{eq:MHD:full}),
the magnetic field equation \eqref{eq:MHD:full:Bdynamics} reduces to
\beq\label{eq:VxB0}
\curl( \vec{u}\times \vec{B})=0 .
\eeq
Beyond the vanishing circulations \eqref{eq:CL:topol:surf:glob},
this conservation law carries important information about
the topology of the magnetic field lines and plasma streamlines.
Where $\vec{u}\times\vec{B}$ can be written as $\grad \psi$,
both families of lines are tangent to the regular level sets $\psi=\const$,
which are consequently \emph{magnetic surfaces};
when compact and without boundary, their connected components are tori,
with nearby level surfaces forming nested families \cite{alexa1935,kruskal1958equilibrium}.

\medskip\noindent\textbf{Example 3.}
Irrotational fluid flows, $\curl\vec{u}=0$,
as well as electrostatic and magnetostatic Maxwell fields,
$\curl\vec{E}=\curl\vec{B}=0$,
yield vanishing circulations of $\vec{u}$, $\vec{E}$, $\vec{B}$
around the boundary of any surface within the solution domain.

Third, a \emph{spatial gradient conservation law} is a vanishing gradient expression
\beq\label{eq:CL:topol:line:loc}
\grad\,\ConsLawFlux=0 ,
\eeq
stating that the scalar flux $\ConsLawFlux[\depvars]$
takes the same value at any two points of the (connected) solution domain;
globally, $\ConsLawFlux\big|_{\p\C}=0$ for any curve $\C$.

\medskip\noindent\textbf{Example 4.}
Bernoulli's principle for steady, inviscid, irrotational, constant-density fluid flows:
setting $\vec{u}_t=0$, $\nu=0$, and $\curl\vec{u}=0$
in the velocity equation \eqref{eq:mom:vxcuv} yields
\beq\label{eq:Bern}
\grad \Big(p +\tfrac{1}{2}|\vec{u}|^2\Big) =0 ,
\eeq
so that the modified pressure \eqref{eq:modified.pressure}
is uniform throughout the flow domain.

Topological conservation laws of all three types
can be sought systematically by an adaptation of the multiplier method
presented in Section~\ref{sec:multipliers},
in which the time-derivative terms are absent;
see \cite{Anc.Che.2020} for details.

\section{Conservation laws on moving domains}\label{sec:CL:moving}

\begin{figure}[htbp]
\begin{subfigure}{0.32\textwidth}
\centering\includegraphics[width=0.98\textwidth]{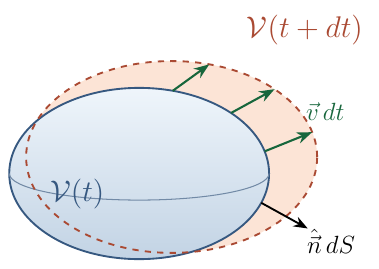}
\caption{}\label{fig:volmov}
\end{subfigure}
\hfill
\begin{subfigure}{0.32\textwidth}
\centering\includegraphics[width=0.98\textwidth]{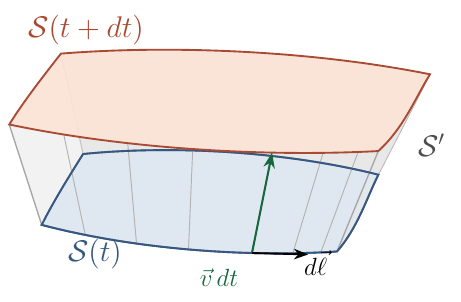}
\caption{}\label{fig:surfmov}
\end{subfigure}
\hfill
\begin{subfigure}{0.32\textwidth}
\centering\includegraphics[width=0.98\textwidth]{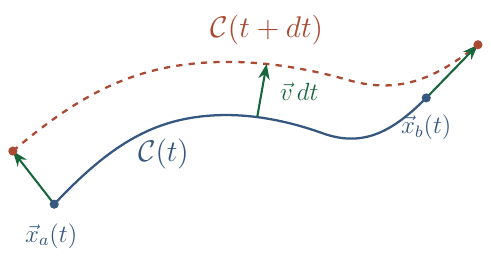}
\caption{}\label{fig:circmov}
\end{subfigure}
\caption{Conservation law computations on moving domains: (\subref{fig:volmov}) a moving volume $\V(t)$; (\subref{fig:surfmov}) a moving surface $\S(t)$ with the swept strip $\S'$; (\subref{fig:circmov}) a moving curve $\C(t)$ with transported endpoints.}
\label{fig:moving}
\end{figure}

The counterparts of divergence-type, curl-type, and gradient-type conservation laws
for moving domains will now be formulated in both global and local forms.
This essentially involves replacing the static domain by
a domain that is transported with respect to a velocity field (which may be time dependent)
and replacing the static flux by a moving flux that takes into account the motion of the domain boundary
(Figure~\ref{fig:moving}). 

\subsection{Moving volumetric conservation laws}\label{sec:CL:moving:vol}

Consider a material domain $\mathcal{V}(t)$ made of labelled particles
that move with the flow induced by a velocity field $\VelVec(\vec{x},t)$.
In particular, their position is given by $\vec{x}(t)$ satisfying
\beq\label{eq:material:particles:traj}
\dfrac{d \vec{x}(t)}{dt} = \VelVec(\vec{x}(t),t),\qquad \vec{x}(0)=\vec{x}_0.
\eeq

A \emph{moving global volumetric conservation law} of a dynamical PDE system 
is a balance equation 
\beq\label{eq:CL:moving:div:global}
\frac{d }{d t} \, \int_{\mathcal{V}(t)} \ConsLawDensity  \;dV = - \int_{\p\mathcal{V}(t)}\, \widetilde{\vec{\ConsLawFlux}} \cdot \dS
\eeq
which holds for all solutions $\depvars(t, \vec{x})$ of the PDE system. 
It has the physical meaning that the rate of change of the total amount of density $\ConsLawDensity$ in the moving domain $\mathcal{V}(t)$ 
is balanced by a flux through the moving boundary $\p\mathcal{V}(t)$.

To establish a local form of the global conservation law \eqref{eq:CL:moving:div:global},
observe that the left hand side of equation \eqref{eq:CL:moving:div:global}
is comprised by two contributions:
the rate of change of the density $\ConsLawDensity(\vec{x},t)$ with time
in the region $\mathcal{V}(t)$, 
and the increment/decrement in the amount contained within the additional region
swept out by the boundary surface $\p\mathcal{V}(t)$
in its motion from time $t$ to $t+dt$
as illustrated in Figure~\ref{fig:volmov}.
The volume of this additional region is $\VelVec\cdot \dS\, dt$. 
As a result, the total rate of change is given by 
\beq\label{eq:Leibniz.vol}
\frac{d }{d t} \, \mathop{\int}_{\mathcal{V}(t)} \ConsLawDensity  \;dV
= \mathop{\int}_{\mathcal{V}(t)} \partial_t \ConsLawDensity  \;dV
+ \mathop{\int}_{\p\mathcal{V}(t)} \ConsLawDensity\,  \VelVec\cdot \dS
= \mathop{\int}_{\mathcal{V}(t)} \left( \partial_t\ConsLawDensity + \div (\ConsLawDensity\, \VelVec) \right)\;dV 
\eeq
after use of the divergence theorem. 
Then, using the right hand side of equation \eqref{eq:CL:moving:div:global}
and again applying the divergence theorem, 
it follows that
\beq
\mathop{\int}_{\mathcal{V}(t)} \left( \partial_t\ConsLawDensity + \div (\ConsLawDensity\, \VelVec) \right)\,dV
=
-\mathop{\int}_{\mathcal{V}(t)} \div\widetilde{\vec{\ConsLawFlux}} \,dV . 
\eeq
(A full proof of this formula is given in Appendix~\ref{sec:proof}.)
This will hold for all moving domains $\mathcal{V}(t)$ if and only if
\beq\label{eq:CL:moving:div:local}
\partial_t \ConsLawDensity + \div\,(\widetilde{\vec{\ConsLawFlux}} + \ConsLawDensity\, \VelVec)=0
\eeq
which constitutes a local divergence-type conservation law
with the flux containing an extra contribution proportional to the velocity field $\VelVec$.

\medskip\noindent\textbf{Example 1.}
In the gas dynamics system \eqref{eq:gas.dyn},
the local volumetric conservation laws of mass \eqref{eq:gas:CL:mass} and entropy \eqref{eq:gas:CL:entropy}
as they stand describe local moving conservation laws
with vanishing moving fluxes,
where the transport velocity is taken to be the gas velocity, $\VelVec=\vec{u}$. 
The resulting global conservation laws on any domain $\V(t)$ transported by the gas flow
represent conservation of total mass and total entropy in a material volume transported by the gas:
\beq
\mathcal{M}= \int_{\V(t)}  \rho \,\dV =\const,
\qquad
\mathcal{S}= \int_{\V(t)}  \rho\, S \,\dV =\const.
\eeq

\medskip\noindent\textbf{Example 2.}
For gas flow \eqref{eq:gas.dyn}, or fluid flow \eqref{eq:NS3d:constd:consf}
in the inviscid case $\nu=0$,
the local volumetric conservation law of energy \eqref{eq:gas:CL:energy} or \eqref{eq:Eul:CL:KinEn}
has the form \eqref{eq:CL:moving:div:local} of a moving conservation law
with $\VelVec=\vec{u}$, 
whose moving flux is $\widetilde{\vec{\ConsLawFlux}} = p\vec{u}$.
Thus, the total energy in any domain $\V(t)$ transported by the flow
\beq\label{eq:KinEner:movdom}
\mathcal{E}= \int_{\V(t)}  \big( \tfrac{1}{2} \rho |\vec{u}|^2 + \rho\, e \big) \,\dV
\quad\text{ or }\quad
\mathcal{E}= \int_{\V(t)}  \tfrac{1}{2} |\vec{u}|^2 \,\dV
\eeq
satisfies 
\beq
\dfrac{d}{dt} \mathcal{E} =  - \int_{\p\V(t)}\, p\vec{u} \cdot \dS
\eeq
where $p\vec{u}$ is the energy flux through the moving boundary $\p\V(t)$.

\medskip\noindent\textbf{Example 3.}
Local conservation law of helicity \eqref{eq:Eul:CL:heli}
for incompressible, constant-density, inviscid fluid flow \eqref{eq:NS3d:constd:consf}
(where $\nu=0$)
can be rewritten by expanding the double cross-product term in the flux,
which gives
\beq\label{eq:Eul:CL:heli:mov}
h_t + \nabla\cdot\big( \vec{u}\times \nabla (\tfrac{1}{2}|\vec{u}|^2 +p) - \vec{\omega} |\vec{u}|^2  + h \vec{u}\big)=0
\eeq
This has the form \eqref{eq:CL:moving:div:local} of a moving volumetric conservation law with $\VelVec=\vec{u}$, where the moving flux is given by 
$\widetilde{\vec{\ConsLawFlux}} = \vec{u}\times \nabla (\tfrac{1}{2}|\vec{u}|^2 +p) - \vec{\omega} |\vec{u}|^2$.
The corresponding global moving conservation law shows that the moving helicity
\beq\label{eq:Heli:movdom}
\mathcal{H}= \int_{\V(t)}  h \,\dV
\eeq
satisfies 
\beq\label{eq:consHeli:movdomain}
\dfrac{d}{dt} \mathcal{H} =
-\int_{\p\V(t)}\, (\vec{u}\times \nabla(\tfrac{1}{2}|\vec{u}|^2 +p) - \vec{\omega} |\vec{u}|^2) \cdot \dS 
\eeq
on any domain $\V(t)$ transported by the fluid. 

In general,
if the velocity field is zero, $\VelVec\equiv0$,
then a moving conservation law \eqref{eq:CL:moving:div:global} or \eqref{eq:CL:moving:div:local}
coincides with an ordinary volumetric conservation law \eqref{3Dglobalvolconslaw}
or \eqref{3Dvolconslaw},
where the moving flux reduces to the static flux
$\widetilde{\vec{\ConsLawFlux}} = \vec{\ConsLawFlux}$.
Conversely, 
any local volumetric conservation law \eqref{3Dvolconslaw}
holding on a static domain $\V$ can be expressed as an equivalent conservation law
on a moving domain $\V(t)$ in a velocity field $\VelVec(t,\vec{x})\not\equiv0$. 
A comparison of equations \eqref{3Dvolconslaw} and \eqref{eq:CL:moving:div:local}
shows that the moving flux is related to the static flux by
\beq\label{eq:CLflux:mov.vol}
\widetilde{\vec{\ConsLawFlux}} = \vec{\ConsLawFlux} - \ConsLawDensity\,\VelVec . 
\eeq

\medskip\noindent\textbf{Example 4.}
For the Keller--Segel model \eqref{eq:KSiegel},
the bacterial-mass conservation law \eqref{eq:KSiegel:VolCL}
has the flux
$\vec{\ConsLawFlux}_\rho=-\nu_\rho\grad\rho+k\rho\grad c$.
For positive $\rho$, define the cell-flux velocity
\beq
\vec{v}_{\rm cell}=\frac{\vec{\ConsLawFlux}_\rho}{\rho}
=-\nu_\rho\grad\log\rho+k\grad c .
\eeq
On any domain $\V(t)$ transported by $\vec{v}_{\rm cell}$,
the moving flux vanishes, and hence the total bacterial population within the domain
is conserved:
\beq
\int_{\V(t)}\rho\,dV=\const. 
\eeq

\medskip\noindent\textbf{Example 5.}
In nonlinear acoustics \eqref{eq:acoustic}, 
the local conservation of helicity 
\beq
h_t + \tfrac{1}{\rho_0} \div \big((p -\nu\,\div\vec{u}) \curl \vec{u}\big) =0
\eeq
holding in the case of zero shear viscosity, $\mu=0$,
yields an equivalent moving conservation law 
\beq
\dfrac{d}{dt} \mathcal{H} =
- \int_{\p\V(t)}\, \Big( \frac{1}{\rho_0}(p -\nu\,\div\vec{u}) \curl \vec{u} -h\vec{u} \Big)\cdot \dS 
\eeq
for helicity \eqref{eq:Heli:movdom} on any domain $\V(t)$ transported by the gas flow.

\medskip\noindent\textbf{Example 6.}
Maxwell's equations \eqref{eq:Max:vac} involve no physically distinguished velocity field.
Let $\VelVec(\vec{x},t)$ therefore be an \emph{arbitrary} velocity field,
which transports the integration domain while the fields and time
are evaluated in the original inertial frame.
By the flux relation \eqref{eq:CLflux:mov.vol},
the energy conservation law \eqref{eq:Max:CL:energy} yields the moving conservation law
\beq\label{eq:Max:CL:moving:energy}
\dfrac{d}{dt} \int_{\V(t)} \tfrac{1}{2}\big(|\vec{E}|^2+|\vec{B}|^2\big)\, \dV
= - \int_{\p\V(t)} \Big( c\,\vec{E}\times\vec{B}
- \tfrac{1}{2}\big(|\vec{E}|^2+|\vec{B}|^2\big)\, \VelVec \Big)\cdot \dS
\eeq
on domains $\V(t)$ transported by $\VelVec$:
the moving flux is the Poynting flux
diminished by the advective flux of the energy density
carried through the moving boundary.

\subsection{Moving surface-flux conservation laws}\label{sec:CL:moving:surf}

Consider a surface $\mathcal{S}(t)$ comprised of material points moving with the flow defined by a velocity field $\vec{v}(\vec{x},t)$.
A \emph{moving global surface-flux conservation law} of a dynamical PDE system 
is a balance equation 
\begin{equation}\label{eq:CL:moving:surf:global}
  \frac{d}{dt} \int_{\mathcal{S}(t)} \vec{\ConsLawDensity} \cdot \dS
  = -\oint_{\p \mathcal{S}(t)} \widetilde{\vec{\ConsLawFlux}} \cdot \,d\vec{\ell}
\end{equation}
with the physical meaning that the rate of change of total flux of the vector field $\vec{\ConsLawDensity}[\depvars]$ through the moving surface ${\mathcal{S}(t)}$
is compensated by the circulation of $\widetilde{\vec{\ConsLawFlux}}$
around the moving boundary curve $\p \mathcal{S}(t)$.

The left-hand side of equation \eqref{eq:CL:moving:surf:global} is a sum of two contributions, 
one being 
the rate of change of the vector density $\vec{\ConsLawDensity}(\vec{x},t)$ with time
in the surface,
and the other being
the increment/decrement in the surface flux due to the motion of the surface $\mathcal{S}$ from $t$ to $t+dt$
which can be decomposed into normal and tangential components.
As seen in Figure~\ref{fig:surfmov},
the normal motion sweeps out a volume with $dV = (\VelVec\,dt)\cdot \dS$,
which contributes a term
$\left(\int_{\mathcal{S}(t)} (\div\vec{\ConsLawDensity}) \VelVec\cdot\dS\right) dt$,
while the motion of the boundary sweeps out a lateral strip with $\dS = (\VelVec\,dt)\times d\vec{\ell}$,
which contributes a term 
$\left( \oint_{\p \mathcal{S}(t)} (\vec{\ConsLawDensity}\times \VelVec ) \cdot d\vec{\ell}\right) dt$.
This leads to the transport identity
\beq\label{eq:Leibniz.surf.gen}
\dfrac{d}{dt} \int_{\mathcal{S}(t)} \vec{\ConsLawDensity} \cdot \dS
= \int_{\mathcal{S}(t)} \big( \partial_t \vec{\ConsLawDensity} + (\div \vec{\ConsLawDensity})\, \VelVec \big) \cdot \dS
+ \oint_{\p \mathcal{S}(t)} (\vec{\ConsLawDensity} \times \VelVec) \cdot d\vec{\ell} \,,
\eeq
(More details of the derivation of this formula are provided in Appendix~\ref{sec:proof}.)
Combining the boundary integral and the right-hand side of equation \eqref{eq:CL:moving:surf:global},
followed by using Stokes' theorem,
then gives
\beq
\int_{\mathcal{S}(t)} \big( \partial_t \vec{\ConsLawDensity} + (\div \vec{\ConsLawDensity})\, \VelVec + \curl\,( \widetilde{\vec{\ConsLawFlux}} + \vec{\ConsLawDensity}\times \VelVec ) \big) \cdot \dS = 0. 
\eeq
This surface integral holds for all moving surfaces
if and only if
\beq\label{eq:CL:moving:curl:local}
\partial_t \vec{\ConsLawDensity}
+(\div\vec{\ConsLawDensity})\VelVec
+ \curl\,( \widetilde{\vec{\ConsLawFlux}} + \vec{\ConsLawDensity}\times \VelVec )=0 ,
\eeq
which yields the local form of the balance equation
for an arbitrary vector density $\vec{\ConsLawDensity}$. 

Taking the divergence of the preceding local balance equation
gives a volumetric transport equation
\beq\label{eq:surfdens.transport}
\partial_t q+\div(q\VelVec)=0,
\quad
q=\div\vec{\ConsLawDensity} . 
\eeq
Consequently, if 
\beq\label{eq:surfdens.divfree}
\div \vec{\ConsLawDensity} = 0
\eeq
holds at any instant, then it is preserved for all time under the flow. 
For such divergence-free vector densities, 
the transport identity \eqref{eq:Leibniz.surf.gen} reduces to the form 
\beq\label{eq:Leibniz.surf}
\dfrac{d}{dt} \int_{\mathcal{S}(t)} \vec{\ConsLawDensity} \cdot \dS
= \int_{\mathcal{S}(t)} \partial_t \vec{\ConsLawDensity} \cdot \dS
+ \oint_{\p \mathcal{S}(t)} (\vec{\ConsLawDensity} \times \VelVec) \cdot d\vec{\ell},
\eeq
and the local balance equation \eqref{eq:CL:moving:curl:local}
correspondingly reduces to a curl-type local conservation law
\beq\label{eq:CL:moving:curl:divfree:local}
\partial_t \vec{\ConsLawDensity}
+ \curl\,( \widetilde{\vec{\ConsLawFlux}} + \vec{\ConsLawDensity}\times \VelVec )=0.
\eeq
It follows that every local vector conservation law of this form \eqref{eq:CL:moving:curl:divfree:local}
for a given PDE system
corresponds to a moving surface-flux conservation law \eqref{eq:CL:moving:surf:global}, 
where the surface flux of $\vec{\ConsLawDensity}$ through a moving surface $\mathcal{S}(t)$
is balanced purely by a boundary circulation integral.
In particular, if $\mathcal{S}(t)$ is a closed surface with no boundary,
for instance, a topological sphere or a torus transported by a velocity field $\VelVec(\vec{x},t)$,
then the total surface flux is conserved, 
\beq
 \int_{\mathcal{S}(t)} \vec{\ConsLawDensity} \cdot \dS  = \const.
 \eeq
The divergence-free condition \eqref{eq:surfdens.divfree} 
holds in all of the examples below.

\medskip\noindent\textbf{Example 1.}
For ideal MHD equations with zero resistivity, $\eta=0$,
the magnetic field evolution equation \eqref{eq:MHD:full:Bdynamics} reduces to 
\beq\label{eq:MHD:Bt:ideal}
\vec{B}_t + \curl( \vec{B}\times \vec{u}) =0,
\eeq
which is exactly in the form \eqref{eq:CL:moving:curl:local}
with zero moving flux $\widetilde{\vec{\ConsLawFlux}}=0$.
As a result, for any surface $\mathcal{S}(t)$ transported by the MHD flow,
the total magnetic flux is conserved 
\beq
 \mathcal{F}_B = \int_{\mathcal{S}(t)} \vec{B} \cdot \dS  = \const,
\eeq
which has the physical meaning that the magnetic field lines are \emph{frozen-in} the moving plasma.

The same holds for the vorticity equation \eqref{eq:NS3d:constd:consf:vort:momom}
in incompressible constant-density fluid flow in the inviscid case $\nu=0$, 
resulting in frozen-in vorticity lines, with the vorticity flux satisfying
\beq
 \mathcal{F}_{\vec{\omega}} = \int_{\mathcal{S}(t)} \vec{\omega} \cdot \dS  = \const
\eeq
for an arbitrary surface $\mathcal{S}(t)$ transported by the fluid flow.

\medskip\noindent\textbf{Example 2.}
In non-ideal MHD \eqref{eq:MHD:full},
the magnetic field evolution equation \eqref{eq:MHD:full:Bdynamics} can be written as
\beq\label{eq:idealMHD.Beqn}
\vec{B}_t + \curl( \vec{B}\times \vec{u} + \eta \curl \vec{B}) =0,
\eeq
by using $\Delta \vec{B} = -\curl\,(\curl \vec{B})$
which holds due to the vector Laplacian identity \eqref{eq:vector:Lapl}
combined with $\div\vec{B}=0$. 
Note that the quantity $\tfrac{1}{\mu} \curl \vec{B} =\vec{J}$ 
is the electric current density vector in the plasma.
Thus, equation \eqref{eq:idealMHD.Beqn} has the form
of a local moving surface-flux conservation law \eqref{eq:CL:moving:curl:local}
with $\widetilde{\vec{\ConsLawFlux}} =\eta\mu\,\vec{J}$.
Its global form gives conservation of magnetic flux 
\beq \label{eq:CL:moving:surf:global:Bgen}
\dfrac{d}{dt} \int_{\mathcal{S}(t)} \vec{B} \cdot d\vec{S}
= -\eta\mu \oint_{\p \mathcal{S}(t)} \vec{J} \cdot \,d\vec{\ell}
\eeq
on moving surfaces $\mathcal{S}(t)$ transported by the plasma flow.

\medskip\noindent\textbf{Example 3.}
Similarly, the vorticity dynamics equation \eqref{eq:NS3d:constd:consf:vort:momom}
in the incompressible constant-density fluid flow
yields the local surface-flux conservation law
\beq
\vec{\omega}_t + \curl( \vec{\omega}\times \vec{u} + \nu \curl \vec{\omega}) =0
\eeq
This describes global conservation of the vorticity flux 
\beq \label{eq:CL:moving:surf:global:Omgen}
\dfrac{d}{dt} \int_{\mathcal{S}(t)} \vec{\omega} \cdot d\vec{S}
= -\nu \oint_{\p \mathcal{S}(t)} \curl \vec{\omega} \cdot \,d\vec{\ell} 
\eeq
on moving surfaces $\mathcal{S}(t)$ transported by the fluid.

As in the volumetric case,
a zero velocity field $\VelVec\equiv0$ reduces the moving conservation law
\eqref{eq:CL:moving:surf:global}, \eqref{eq:CL:moving:curl:local}
to its static counterpart \eqref{3Dglobalsurfconslaw}, \eqref{3Dsurfconslaw};
conversely, any local surface-flux conservation law
with a divergence-free density \eqref{eq:surfdens.divfree}
can equivalently be formulated on moving surfaces,
with the moving and static fluxes related by
\beq\label{eq:CLflux:mov.surf}
\widetilde{\vec{\ConsLawFlux}} = \vec{\ConsLawFlux} - \vec{\ConsLawDensity}\times \VelVec .
\eeq

\medskip\noindent\textbf{Example 4.}
In nonlinear acoustics \eqref{eq:acoustic},
the vorticity equation \eqref{eq:acoustic.vort} can be written in the form 
\beq\label{eq:acoustic.vort:surf}
 \vec{\omega}_t + \frac{\mu}{\rho_0}\, \curl (\curl \vec{\omega}) =0
\eeq
after using the vector Laplacian identity \eqref{eq:vector:Lapl}
and $\div\vec{\omega}=0$. 
This local conservation law is equivalent to a moving surface-flux conservation law \eqref{eq:CL:moving:curl:local}
with $\VelVec=\vec{u}$ and $\widetilde{\vec{\ConsLawFlux}} = \dfrac{\mu}{\rho_0}\, \curl \vec{\omega} - \vec{\omega}\times \vec{u}$.
The corresponding global conservation law is given by 
\beq \label{eq:CL:moving:surf:global:acoustic}
\dfrac{d}{dt} \int_{\mathcal{S}(t)} \vec{\omega} \cdot d\vec{S}
= -\oint_{\p \mathcal{S}(t)} \Big( \frac{\mu}{\rho_0}\, \curl \vec{\omega} - \vec{\omega}\times \vec{u} \Big)\cdot d\vec{\ell} 
\eeq
which describes conservation of vorticity flux 
on any moving surface $\mathcal{S}(t)$ transported by the gas flow.

\medskip\noindent\textbf{Example 5.}
For Maxwell's equations \eqref{eq:Max:vac},
again take $\VelVec(\vec{x},t)$ to be an arbitrary velocity field
transporting the integration surface in the original inertial frame.
Both field equations \eqref{eq:Max:Edyn}--\eqref{eq:Max:Bdyn}
give local moving surface-flux conservation laws,
with the divergence-free condition \eqref{eq:surfdens.divfree}
holding exactly by virtue of the constraints \eqref{eq:Max:divE}--\eqref{eq:Max:divB}.
For the magnetic field,
the moving flux is
$\widetilde{\vec{\ConsLawFlux}} = c\vec{E} - \vec{B}\times\VelVec$,
and the resulting global conservation law
\beq\label{eq:Max:CL:moving:Faraday}
\dfrac{d}{dt} \int_{\mathcal{S}(t)} \vec{B} \cdot \dS
= -\oint_{\p \mathcal{S}(t)} \big( c\vec{E} + \VelVec\times\vec{B} \big) \cdot d\vec{\ell}
\eeq
is Faraday's law of induction for a moving circuit.
In the units used here, the electromotive force is
$\oint_{\p\mathcal{S}(t)}(\vec{E}+c^{-1}\VelVec\times\vec{B})\cdot d\vec{\ell}$;
the second term accounts for motional induction.
Likewise, for the electric field,
$\widetilde{\vec{\ConsLawFlux}} = -c\vec{B} - \vec{E}\times\VelVec$
yields
\beq\label{eq:Max:CL:moving:Ampere}
\dfrac{d}{dt} \int_{\mathcal{S}(t)} \vec{E} \cdot \dS
= \oint_{\p \mathcal{S}(t)} \big( c\vec{B} - \VelVec\times\vec{E} \big) \cdot\dl, 
\eeq
the source-free Amp\`ere--Maxwell law for a moving circuit.

\subsection{Moving circulatory conservation laws}\label{sec:CL:moving:circ}

Finally, consider a curve $\C(t)$ comprised of material points moving with the flow defined by a velocity field $\vec{v}(\vec{x},t)$.
A \emph{moving global circulatory conservation law} of a dynamical PDE system 
is a balance equation 
\begin{equation}\label{eq:CL:moving:curv:global}
  \frac{d}{dt} \int_{\C(t)} \vec{\ConsLawDensity} \cdot \dl
  = -\widetilde{\ConsLawFlux}\Big|_{\p\C(t)}
\end{equation}
which has the physical meaning that the rate of change of
total circulation of the vector field $\vec{\ConsLawDensity}[\depvars]$ around the moving curve $\C(t)$
is compensated by the net outward flow through the endpoints $(\vec{x}_a(t),\vec{x}_b(t))$ of the curve.

Similarly to a moving surface conservation law, 
the left-hand side of equation \eqref{eq:CL:moving:curv:global}
is a sum of contributions from the rate of change of the circulation density with time
and the motion of the curve from $t$ to $t+dt$
as shown in Figure~\ref{fig:circmov}.
Applying Stokes' theorem to the surface swept out by the curve gives a curl contribution
and contributions from both endpoint trajectories.
The latter involve the full endpoint velocities, including their normal components.
This leads to the transport identity
\beq\label{eq:Leibniz.curv.gen}
\dfrac{d}{dt} \int_{\C(t)} \vec{\ConsLawDensity}   \cdot d\vec{\ell}
=  \int_{\C(t)} \big( \partial_t \vec{\ConsLawDensity} + (\curl \vec{\ConsLawDensity})\times\VelVec \big) \cdot d\vec{\ell}
+\big(\vec{\ConsLawDensity} \cdot \VelVec\big) \Big|_{\vec{x}_a(t)}^{\vec{x}_b(t)}
\eeq
(More details of the derivation of this formula are provided in Appendix~\ref{sec:proof}.)
The line integral combined with the right-hand side of equation \eqref{eq:CL:moving:curv:global}
yields
\beq\label{eq:curv.transport.gen}
\int_{\C(t)} \big( \partial_t \vec{\ConsLawDensity} + (\curl \vec{\ConsLawDensity})\times\VelVec \big) \cdot d\vec{\ell}
+\big(\vec{\ConsLawDensity} \cdot \VelVec + \widetilde{\ConsLawFlux} \big) \Big|_{\vec{x}_a(t)}^{\vec{x}_b(t)}
= \int_{\C(t)} \big( \partial_t \vec{\ConsLawDensity} + (\curl \vec{\ConsLawDensity})\times\VelVec + \grad(\vec{\ConsLawDensity} \cdot \VelVec + \widetilde{\ConsLawFlux}) \big)
\cdot d\vec{\ell}
=0
\eeq
after use of the fundamental theorem of calculus for line integrals.
This equation shows that the global balance law \eqref{eq:CL:moving:curv:global}
holds for all moving curves if and only if
\begin{equation}\label{3Dcurvconslaw:mov}
\partial_t \vec{\ConsLawDensity}
+(\curl\vec{\ConsLawDensity})\times\VelVec
+ \grad\,(\widetilde{\ConsLawFlux} + \vec{\ConsLawDensity}\cdot \VelVec) =0.
\end{equation}
which yields the local form of the balance law
for an arbitrary circulation density $\vec{\ConsLawDensity}$. 

Taking the curl of the preceding local balance equation
gives a surface-flux transport equation 
\beq\label{eq:curldens.transport}
\partial_t\vec{\omega}+\curl(\vec{\omega}\times\VelVec)=0,
\quad
\vec{\omega}=\curl\vec{\ConsLawDensity}
\eeq
Consequently, irrotationality is preserved whenever it holds at one instant.

A sufficient condition for a direct reduction of the local balance law \eqref{3Dcurvconslaw:mov}
to a gradient-type conservation law consists of 
\beq\label{eq:curvdens.curlfree}
(\curl \vec{\ConsLawDensity})\times\VelVec = 0.
\eeq
This holds, in particular, when $\vec{\ConsLawDensity}$ is irrotational
or when its curl is parallel to $\VelVec$.
Through condition \eqref{eq:curvdens.curlfree},
the global transport identity \eqref{eq:Leibniz.curv.gen} reduces to the form 
\beq\label{eq:Leibniz.curv}
\dfrac{d}{dt} \int_{\C(t)} \vec{\ConsLawDensity} \cdot\dl
= \int_{\C(t)} \partial_t \vec{\ConsLawDensity} \cdot \dl
+\big(\vec{\ConsLawDensity}\cdot\VelVec\big)\Big|_{\vec{x}_a(t)}^{\vec{x}_b(t)}
\eeq
For any closed moving loop, the general balance law
\eqref{eq:CL:moving:curv:global} has no end-point flux, independently of
condition \eqref{eq:curvdens.curlfree}, and hence 
\begin{equation}\label{eq:CL:moving:curv:global.com}
  \int_{\C(t)} \vec{\ConsLawDensity} \cdot \dl =\const.
\end{equation}

As in the cases of moving volume and surface conservation laws, 
a zero velocity field $\VelVec\equiv0$ reduces the moving conservation law
\eqref{eq:CL:moving:curv:global}, \eqref{3Dcurvconslaw:mov}
to its static counterpart \eqref{3Dglobalcurvconslaw}, \eqref{3Dcurvconslaw};
conversely, any local circulatory conservation law
with a density satisfying condition \eqref{eq:curvdens.curlfree}
can equivalently be formulated on moving curves,
with the moving and static end-point fluxes related by
\beq\label{eq:CLflux:mov.curv}
\widetilde{\ConsLawFlux} = \ConsLawFlux - \vec{\ConsLawDensity}\cdot \VelVec .
\eeq

\medskip\noindent\textbf{Example 1.}
For nonlinear acoustics \eqref{eq:acoustic} in the case $\mu=0$, 
the momentum equation \eqref{eq:acoustic.vel} reduces to
a local circulatory conservation law
\beq
\vec{u}_t = - \grad\Big(  \frac{1}{\rho_0} (p  - \nu\,\div \vec{u}) \Big) .
\eeq
Consider irrotational flows, $\vec{\omega}=\curl\vec{u}=0$,
which persist in time when $\mu=0$ by equation \eqref{eq:acoustic.vort},
so that the condition \eqref{eq:curvdens.curlfree} holds.
Then, on curves transported by the gas flow, $\VelVec=\vec{u}$,
this is a moving circulatory conservation law
with the moving end-point flux given by
$\widetilde\ConsLawFlux = \frac{1}{\rho_0} (p  - \nu\,\div \vec{u}) - |\vec{u}|^2$.
The global form of this moving conservation law
\beq
\frac{d}{dt}\int_{\C(t)} \vec{u} \cdot d\vec{\ell}
= - \Big( \frac{1}{\rho_0} (p  - \nu\,\div \vec{u})  - |\vec{u}|^2 \Big)\Big|_{\p\C(t)}
\eeq
describes the rate of change of total circulation
as being equal to minus the net outward flux through the end points of the moving curve.

\medskip\noindent\textbf{Example 2.}
A similar example is given by the case of Beltrami flows,
$\curl\vec{u} = \alpha \vec{u}$,
in incompressible constant-density fluid flow \eqref{eq:NS3d:constd:consf}
with zero viscosity $\nu=0$.
The resulting momentum equation \eqref{eq:mom:vxcuv}
combined with the identity \eqref{eq:importantID:curlcurl}
yields
\beq\label{eq:mom:vxcuv:circ:mov}
\vec{u}_t + \grad \Big(p + \frac{1}{2} |\vec{u}|^2\Big) =0
\eeq
This has the form of a local moving circulatory conservation law \eqref{3Dcurvconslaw:mov}
with $\VelVec=\vec{u}$,
where $\widetilde\ConsLawFlux = p - \frac{1}{2} |\vec{u}|^2$
is the moving end-point flux.
Note that the condition \eqref{eq:curvdens.curlfree} holds automatically for Beltrami flows,
since $(\curl\vec{u})\times\vec{u} = \alpha\, \vec{u}\times\vec{u} = 0$.
In global form, 
\beq
\dfrac{d}{dt} \int_{\C(t)} \vec{u} \cdot d\vec{\ell} = - \Big(p - \frac{1}{2}|\vec{u}|^2 \Big) \Big|_{\p\C(t)}
\eeq
describes the rate of change of total circulation
holding for curves $\C(t)$ transported by the fluid.

\section{Origin and construction of local conservation laws in terms of multipliers}\label{sec:multipliers}

Under the regularity assumptions stated below, all three types of local conservation laws
--- divergence-type, curl-type, and gradient-type, on static and moving domains ---
can be constructed using the multiplier method and its extensions developed here.
The PDE system is assumed to admit a regular solved form, with compatibility conditions included
(see Refs.~\cite{olver1993applications,bluman2010applications,Anc-review}).
This setting encompasses Kovalevskaya and generalized Kovalevskaya systems,
as well as the constrained systems considered below.

A general explanation of how multipliers characterize
local conservation laws will now be given in three dimensions.
Consider a general dynamical PDE system of $N\geq 1$ equations 
\beq\label{sys}
  G[\depvars]=(G^1[\depvars],\ldots,G^N[\depvars])=0
\eeq
with $m\geq1$ dependent variables
\beq
\depvar = (\depvar^1,\ldots,\depvar^m)
\eeq
which are functions of time $t$ and position $\vec{x}=(x^1,x^2,x^3)$.
The spatial coordinates may be any positively oriented, time-independent
curvilinear coordinates (Appendix~\ref{sec:curvilin.coords}).
Write $g_{ij}$ for the metric components and $J=\sqrt{\det g}$ for the volume factor.
Thus, for a vector differential function $\vec F$,
\[
\div\vec F=\frac{1}{J}D_i(JF^i).
\]
Derivatives acting on differential functions are total derivatives;
$D_t$ and $D_i$ will be used when this distinction needs to be explicit. 
Each equation $G^\sigma[\depvars]$ in the system is assumed to be described
by a differential function of $\depvars$.
After a change of variables if necessary,
the system will be required to possess a set of $N$ leading derivatives such that
the PDEs have the form 
\beq\label{sys.solvedform}
\begin{aligned}
  G^1[\depvars] & = \partial_{s^1}^{j_1} \depvar^{\alpha_1} - g^1[\depvars],
  \\
  \vdots & 
  \\
  G^N[\depvars] &= \partial_{s^N}^{j_N} \depvar^{\alpha_N} - g^N[\depvars],
  \end{aligned}
\eeq
where $j_1,\ldots,j_N$ are positive integers;
$\alpha_1,\ldots,\alpha_N$ belong to the set $\{1,\ldots,m\}$;
$s^1,\ldots,s^N$ are specific functions of $t,x^1,x^2,x^3$;
and $g^1[\depvars],\ldots,g^N[\depvars]$ are differential functions of
$(\depvar^1,\ldots,\depvar^m)$
with neither $\partial_{s^1}^{j_1} \depvar^{\alpha_1},\ldots,\partial_{s^N}^{j_N} \depvar^{\alpha_N}$
nor any of their differential consequences
appearing in these functions.
The solved form is assumed to be regular and formally integrable, with all independent compatibility conditions included.

\subsection{Examples of leading derivatives}

A given dynamical PDE system may possess more than one solved form \eqref{sys.solvedform}.
To illustrate the situation, we give choices of leading derivatives for the running examples in Section~\ref{sec:PDEsystems}.
The compact vector expressions below are understood with the constraints and their differential consequences substituted as needed to eliminate leading derivatives from the right-hand sides.

\medskip\noindent\textbf{Maxwell's electromagnetic equations in vacuum}:
For the 6 dependent variables
$\depvars=(E^1,E^2,E^3,B^1,B^2,B^3)$ in Cartesian coordinates,
there are 6 evolution PDEs and two constraint PDEs 
\beq\label{eq:Max:vac:leading}
\vec{E}_t  = c\, \curl \vec{B},
\quad
\vec{B}_t  = - c\, \curl \vec{E},
\quad
E^1_x = -E^2_y - E^3_z,
\quad
B^1_x = -B^2_y - B^3_z,
\eeq
with leading derivatives given in terms of
$s^1=\cdots=s^6=t$, $s^7=s^8=x$, and $j_1=\cdots=j_8=1$.

\medskip\noindent\textbf{Euler equations of compressible gas dynamics}:
For the 5 dependent variables
$\depvars=(\rho,p,u^1,u^2,u^3)$ in Cartesian coordinates,
there are 5 evolution PDEs 
\beq\label{eq:gas.dyn:leading}
\rho_t  = -\div (\rho \vec{u}),
\quad
\vec{u}_t = - (\vec{u}\cdot\nabla)\vec{u} - \frac{1}{\rho}\grad {p},
\quad
p_t = -(\vec{u}\cdot\nabla)p - F(\rho,p) \div \vec{u}
\eeq
with leading derivatives given in terms of
$s^1=\cdots=s^5=t$, and $j_1=\cdots=j_5=1$.

\medskip\noindent\textbf{Keller--Segel model of chemotaxis}:
For the two dependent variables $\depvars = (\rho,c)$,
there are two evolution PDEs
\beq
\rho_t  =\nu_\rho \Delta \rho - k\,\div (\rho \grad c) ,
\quad
c_t =\nu_c \Delta c + g(\rho) - h(c),
\eeq
with leading derivatives given in terms of
$s^1=s^2=t$, and $j_1=j_2=1$.

An alternative set of leading derivatives is $s^1=s^2=x$, and $j_1=j_2=2$, which gives
\beq
\begin{aligned}
\rho_{xx} & = \frac{1}{\nu_\rho}\Big( \rho_t + k\,\grad\rho\cdot \grad c +\frac{k}{\nu_c}\,\rho \big(c_t + h(c) -g(\rho)\big) \Big) -\rho_{yy} -\rho_{zz}, 
\\
c_{xx} & = \frac{1}{\nu_c}\Big( c_t +h(c) - g(\rho) \Big) -c_{yy} -c_{zz}
\end{aligned}
\eeq

\medskip\noindent\textbf{Navier-Stokes equations of incompressible constant-density fluid flow}:
For the four dependent variables $\depvars=(u^1,u^2,u^3,p)$ in Cartesian coordinates,
there is one constraint PDE, three evolution PDEs,
and the pressure Poisson equation,
which arises from the divergence of the velocity equation
combined with the incompressibility constraint:
\beq
\begin{array}{l}
u^1_x = -u^2_y -u^3_z,
\qquad
\vec{u}_t = -(\vec{u}\cdot\nabla)\vec{u} -\grad p + \nu \,\Delta \vec{u},
\\[1.5ex]
p_{xx} = -p_{yy} - p_{zz} - \div\left( (\vec{u}\cdot\nabla)\vec{u} \right)
\end{array}
\eeq
with leading derivatives given in terms of
$s^1=x$, $s^2=s^3=s^4=t$, $s^5=x$, and $j_1=\cdots=j_4=1$, $j_5=2$.

\medskip\noindent\textbf{Nonlinear acoustic equations}:
For the 4 dependent variables $\depvar=(\rho,u^1,u^2,u^3)$ in Cartesian coordinates,
there are 4 evolution equations
\beq
\rho_t = - \rho_0 \div \vec{u} ,
\quad
\vec{u}_t = \dfrac{1}{\rho_0}\left( -\grad p + \mu\, \Delta \vec{u} + \nu\, \grad\div \vec{u} \right),
\eeq
with leading derivatives given in terms of
$s^1=\cdots=s^4=t$, and $j_1=\cdots=j_4=1$.

The wave equation \eqref{eq:acoustic.dens.waveeqn} for $\rho$ is a 2nd-order evolution equation
with leading derivative $s^1=t$ and $j_1 =2$.

\medskip\noindent\textbf{Equations of incompressible magnetohydrodynamics}:
For the 8 dependent variables $\depvars=(\rho,u^1,u^2,u^3,B^1,B^2,B^3,P)$ in Cartesian coordinates,
there are 7 evolution equations, two constraint equations,
and a pressure equation arising from the divergence of the velocity equation
combined with the constraint $\div\vec{u}=0$:
\beq
\begin{array}{l}
\rho_t = - \div (\rho\vec{u}),
\qquad
\vec{u}_t = \vec{u}\times \curl \vec{u} + \dfrac{1}{\mu\rho}\curl \vec{B} \times \vec{B} -\dfrac{1}{\rho}\grad P 
- \tfrac{1}{2} \grad |\vec{u}|^2 + \dfrac{\lambda}{\rho} \Delta \vec{u},
\\[2ex]
\vec{B}_t = \curl(\vec{u} \times \vec{B}) + \eta\, \Delta \vec{B},
\qquad
B^1_x = -B^2_y -B^3_z,
\qquad
u^1_x = -u^2_y -u^3_z,
\\[2ex]
P_{xx} = -P_{yy} - P_{zz} + \dfrac{1}{\rho}\,\grad\rho\cdot\grad P
\\[1.5ex]
\qquad\qquad
+\, \rho\, \div\left( \vec{u}\times \curl \vec{u} - \tfrac{1}{2}\grad|\vec{u}|^2
+ \dfrac{1}{\mu\rho}\curl \vec{B} \times \vec{B} + \dfrac{\lambda}{\rho}\Delta\vec{u} \right)
\end{array}
\eeq
with leading derivatives given in terms of
$s^1=\cdots=s^7=t$, $s^8=s^9=x$, $s^{10}=x$, and $j_1=\cdots=j_9=1$, $j_{10}=2$.

\subsection{Multipliers}

The arguments below are local near a regular point of the prolonged PDE system.
The independent prolonged equations serve as normal coordinates there;
a conservation law is understood to hold on this equation manifold.
The resulting representation of a differential function that vanishes there
is often called Hadamard's lemma:
If a finite-order differential function $f[\depvars]$ vanishes on the regular
prolonged equation manifold of the PDE system,
then locally it can be expressed as a finite linear combination of the PDEs
and their total derivatives, 
\beq\label{eq:hadamard}
f[\depvars]
= \Dop_1(G^1) + \cdots +\Dop_N(G^N)
\eeq
where each $\Dop_1,\cdots,\Dop_N$ is a linear differential operator with respect to $t,\vec{x}$
having coefficients given by differential functions that are non-singular
when evaluated for any solution $\depvars(t,\vec{x})$.
A key corollary of this result is that
\beq\label{eq:multr}
f[\depvars]
= \partial_t \Gamma^t[\depvar;G] + \div \vec{\Gamma}[\depvars;G]
+ Q_1[\depvars] G^1 + \cdots +Q_N[\depvars] G^N
\eeq
holds off solutions, i.e.
when $\depvars(t,\vec{x})$ is an arbitrary function of $t$ and $\vec{x}$, 
which follows by integrating by parts in $Jf$ and then dividing by $J$.
The factors $J$ convert coordinate divergence terms into $\div$ terms.
Here $\Gamma^t[\depvar;G]$ and $\vec{\Gamma}[\depvars;G]$
are differential functions that vanish when $G=0$,
while $Q_1[\depvars], \ldots,Q_N[\depvars]$
are differential functions that are non-singular
when evaluated for solutions $\depvars(t,\vec{x})$ of the given PDE system. 
This set of functions is called a \emph{multiplier}
associated to $f[\depvar]$.

\subsubsection{Local divergence-type conservation laws}

To proceed,
first consider local divergence-type conservation laws \eqref{3Dvolconslaw}
for a given dynamical PDE system. 
Since any conservation law must hold for all solutions of the system, 
it is equivalent to an identity
\beq
\partial_ t\ConsLawDensity + \div \vec{\ConsLawFlux}
= \partial_t \Gamma^t[\depvar;G] + \div \vec{\Gamma}[\depvars;G]
+ Q_1[\depvars] G^1 + \cdots +Q_N[\depvars] G^N
\eeq
holding for arbitrary functions $\depvars(t,\vec{x})$.
This implies
\beq\label{eq:div.conslaw.ident}
\partial_ t\big(\ConsLawDensity - \Gamma^t[\depvar;G] \big)
+ \div\big(\vec{\ConsLawFlux} - \vec{\Gamma}[\depvars;G]\big)
= Q_1[\depvars] G^1 + \cdots +Q_N[\depvars] G^N
\eeq
In particular,
any local volumetric conservation law in a given PDE system is equivalent to
a corresponding local divergence identity \eqref{eq:div.conslaw.ident}.

Moreover, if two conservation laws
given by pairs $(\ConsLawDensity,\ConsLawFlux)_1$ and $(\ConsLawDensity,\ConsLawFlux)_2$ 
lead to identities that have the same right hand side,
then their difference is a trivial conservation law as follows
(see Ref.\cite{olver1993applications,bluman2010applications,Anc-review}). 
The difference in densities will satisfy 
$\ConsLawDensity_1 - \ConsLawDensity_2 = \div \vec{\Upsilon}[\depvars]$
and difference in fluxes will satisfy
$\ConsLawFlux_1 - \ConsLawFlux_2 = \curl\vec{\Lambda}[\depvars]  -\partial_t \vec{\Upsilon}[\depvars]$,
holding for all solutions $\depvars(t,\vec{x})$ of the given PDE system,
with $\vec{\Upsilon}[\depvars]$ and $\vec{\Lambda}[\depvars]$ 
being some specific vector differential functions.

If a PDE system admits a nontrivial differential identity of the form
\[
\Dop_1(G^1)+\cdots+\Dop_N(G^N)\equiv0,
\]
where the $\Dop_\sigma$ are linear differential operators with differential-function
coefficients, then there will exist multipliers that correspond to trivial conservation laws
(see \cite{Anc-review,Anc.Che.2020} for more details).
In contrast, for PDE systems such as evolution equations or Kovalevskaya systems,
which do not possess nontrivial differential identities,
the only multipliers that correspond to trivial conservation laws
are ones that vanish for all solutions of the system.

The simultaneous conditions $\sg E_{\depvar^a}(f)=0$, $a=1,\ldots,m$,
for the Euler operators characterize differential functions that are locally
total coordinate divergences $D_tf^t+D_if^i$.
Since $J$ is time independent, the left-hand side of
\eqref{eq:div.conslaw.ident}, multiplied by $J$, has precisely this form.
Consequently, the multiplier determining equations in the chosen coordinates are
\beq\label{eq:div.conslaw.detsys}
\sg E_{\depvar^a}\!\left(J\sum_{\sigma=1}^NQ_\sigma G^\sigma\right)\equiv0,
\qquad a=1,\ldots,m.
\eeq
They hold for arbitrary functions $\depvars(t,\vec x)$, not just for solutions
of the PDE system. In Cartesian coordinates $J=1$, giving the usual equations.
Each solution of this linear determining system yields a local volumetric
conservation law, whose density and flux can be recovered by standard methods
\cite{olver1993applications,bluman2010applications,Anc-review}.

\subsubsection{Differential operators for surface and curve laws}

For surface and curve laws, integration by parts cannot in general remove
all derivatives of the PDEs while preserving the curl or gradient form.
Therefore, combinations of the equations and their derivatives need to be retained:
\beq\label{eq:operator.characteristic}
\vec b=\sum_{\sigma=1}^N\sum_{|I|\leq r}
\vec R_{\sigma I}[\depvars]D_I G^\sigma.
\eeq
Here $I$ is a spacetime multi-index, $r$ is finite, and the coefficients
$\vec R_{\sigma I}$ are smooth differential functions depending on finitely many
partial derivatives. Such an expression is called a \emph{differential-operator
characteristic}; it vanishes whenever the PDEs and their differential consequences hold.
The transport field $\VelVec$ is given, either as a function of $t,\vec x$
or as a differential function of $\depvars$.

The construction uses the divergence or curl of the unknown density.
For a surface law, set $q=\div\vec{\ConsLawDensity}$;
for a curve law, set $\vec q=\curl\vec{\ConsLawDensity}$.
Taking the divergence or curl of the local balance equation eliminates its flux.
This gives compatibility equations for $q$ or $\vec q$ and the coefficients
in \eqref{eq:operator.characteristic}.
Once these equations are solved, the density and flux can be reconstructed as follows.

\begin{proposition}\label{prop:operator.characteristics}
Under the stated regularity assumptions, every moving surface-flux conservation law
has a representation \eqref{eq:operator.characteristic} with a scalar differential
function $q$ satisfying
\beq\label{eq:operator.surface}
D_tq+\div(q\VelVec)-\div\vec b\equiv0.
\eeq
Conversely, every such pair $(q,\vec b)$ locally determines a moving surface-flux
conservation law with $\div\vec{\ConsLawDensity}=q$.
For a moving circulatory conservation law, the corresponding determining identities are
\beq\label{eq:operator.curve}
\div\vec q\equiv0,
\qquad D_t\vec q+\curl(\vec q\times\VelVec)-\curl\vec b\equiv0.
\eeq
Every moving circulatory law has such a representation, and conversely every
pair satisfying these identities locally determines a law with
$\curl\vec{\ConsLawDensity}=\vec q$.
The reconstructed densities and fluxes are differential functions of finite order.
Setting $\VelVec=0$ gives the static cases.
\end{proposition}

\begin{proof}
Necessity follows by applying \eqref{eq:hadamard} to each component of
\eqref{eq:CL:moving:curl:local} or \eqref{3Dcurvconslaw:mov}, giving the
right-hand side $\vec b$. Taking the divergence or curl then gives
the determining identities.

For the converse, we use the following vector-calculus consequences of
exactness of the differential-function complex
\cite[Theorem~5.80, pp.~352, 369--373]{olver1993applications}.
An identity
\[
D_tq-\div\vec a\equiv0
\]
admits simultaneous potentials $\vec{\ConsLawDensity}$ and $\vec H$ such that
\beq\label{eq:operator.surface.potentials}
q=\div\vec{\ConsLawDensity},
\qquad
\vec a=D_t\vec{\ConsLawDensity}+\curl\vec H.
\eeq
Likewise, the identities
\[
\div\vec q\equiv0,
\qquad D_t\vec q-\curl\vec a\equiv0
\]
admit simultaneous potentials $\vec{\ConsLawDensity}$ and $H$ such that
\beq\label{eq:operator.curve.potentials}
\vec q=\curl\vec{\ConsLawDensity},
\qquad
\vec a=D_t\vec{\ConsLawDensity}+\grad H.
\eeq
These statements apply locally on sufficiently small star-shaped coordinate
neighborhoods of the unrestricted jet bundle, where the field variables and their partial derivatives
are independent coordinates. They are the vector potential statements of the
algebraic Poincare lemma in the four independent variables $t,x^1,x^2,x^3$.
The lemma supplies potentials depending on finitely many derivatives, rather
than spatial or temporal integrals of a particular solution.
The time-independent metric and volume factor convert these potential equations
between Cartesian and curvilinear coordinates.

For \eqref{eq:operator.surface}, take $\vec a=\vec b-q\VelVec$ in
\eqref{eq:operator.surface.potentials} and define
$\widetilde{\vec{\ConsLawFlux}}=\vec H-\vec{\ConsLawDensity}\times\VelVec$.
For \eqref{eq:operator.curve}, take
$\vec a=\vec b-\vec q\times\VelVec$ in \eqref{eq:operator.curve.potentials}
and define $\widetilde{\ConsLawFlux}=H-\vec{\ConsLawDensity}\cdot\VelVec$.
The reconstructed equations are, respectively,
\beq\label{eq:operator.reconstruction}
\begin{aligned}
D_t\vec{\ConsLawDensity}
+(\div\vec{\ConsLawDensity})\VelVec
+\curl(\widetilde{\vec{\ConsLawFlux}}+\vec{\ConsLawDensity}\times\VelVec)&=\vec b,\\
D_t\vec{\ConsLawDensity}
+(\curl\vec{\ConsLawDensity})\times\VelVec
+\grad(\widetilde{\ConsLawFlux}+\vec{\ConsLawDensity}\cdot\VelVec)&=\vec b.
\end{aligned}
\eeq
Since $\vec b=0$ on solutions, these are the required local moving conservation laws.
\end{proof}

Here $q$ and $\vec q$ are unknown differential functions in the determining
identities, not extra constraints imposed on solutions. In particular, their
values may differ between solutions. The determining equations are linear in
these functions and in the coefficients $\vec R_{\sigma I}$ for a given
transport field. A computation begins by choosing finite differential-order
and dependence classes for these unknowns. Differential identities among the
PDEs can make the representation nonunique.

\subsubsection{Local curl-type conservation laws}

For static curl-type laws, a useful special ansatz in
Proposition~\ref{prop:operator.characteristics} is
$q=\xi(\vec x)+\sum_{\sigma}Q_\sigma G^\sigma$ and
$\vec b=\sum_{\sigma}\vec Q_\sigma G^\sigma$,
where $\xi$ is a prescribed field-independent spatial function.
After the compatibility condition below is imposed, the reconstructed identities are
\beq\label{eq.surf.conslaw.ident1}
\div\vec{\ConsLawDensity}
 =\xi(\vec x)+Q_1G^1+\cdots+Q_NG^N
\eeq
and
\beq\label{eq.surf.conslaw.ident2}
\partial_t\vec{\ConsLawDensity}+\curl\vec{\ConsLawFlux}
 =\vec Q_1G^1+\cdots+\vec Q_NG^N.
\eeq
These identities hold for arbitrary functions $\depvars(t,\vec x)$.
The scalar and vector functions $Q_\sigma$ and $\vec Q_\sigma$
are algebraic multipliers: no derivatives act on $G^\sigma$ in these expressions.
This ansatz is sufficient but is not asserted to represent every curl-type law.
In particular, time independence of $\div\vec{\ConsLawDensity}$ on solutions
does not require its value to be the same for different solutions.
For example, $u_t=0$ admits the surface-flux law
$\vec{\ConsLawDensity}=(u,0,0)$, $\vec{\ConsLawFlux}=0$.
Here $\div\vec{\ConsLawDensity}=u_x$ is time independent, but varies with the
arbitrary stationary solution $u(\vec x)$.

The complete compatibility condition for this ansatz is
\beq\label{eq:curl.conslaw.deteq}
\partial_t(Q_1 G^1 +\cdots + Q_NG^N) -\div(\vec{Q}_1 G^1 + \cdots +\vec{Q}_N G^N)
\equiv 0
\eeq
which must hold for arbitrary functions $\depvars(t,\vec{x})$.
This linear determining equation will split with respect to derivatives of $\depvars$
that do not appear in the multiplier functions.
The resulting equations imply the necessary Euler tests
\beq\label{eq:curl.conslaw.detsys}
\begin{aligned}
&\widehat{\sg E}_{\depvar^a}\!\left(J\sum_\sigma Q_\sigma G^\sigma\right)\equiv0,\\
&\sg E_{\depvar^a}\!\left(J\sum_\sigma Q_\sigma^{\,i}G^\sigma\right)\equiv0,
\qquad a=1,\ldots,m,\quad i=1,2,3.
\end{aligned}
\eeq
The factors $J$ follow from the component identities
$J\div\vec{\ConsLawDensity}=D_i(J\ConsLawDensity^i)$ and
$J(D_t\vec{\ConsLawDensity}+\curl\vec{\ConsLawFlux})^i
=D_t(J\ConsLawDensity^i)+\epsilon^{ijk}D_j\ConsLawFlux_k$,
where $\ConsLawFlux_k=g_{kl}\ConsLawFlux^l$.
Here $\widehat{\sg E}_{\depvar^a}$ denotes the family of spatial Euler operators
\beq\label{eq:spatial.euler.family}
\widehat{\sg E}_{D_t^k\depvar^a}
 =\sum_K(-D_x)_K\frac{\partial}{\partial(D_KD_t^k\depvar^a)},
\qquad k=0,1,\ldots,
\eeq
where $K$ ranges over spatial multi-indices and only finitely many terms
act nontrivially on a given differential function.
All temporal jets $D_t^k\depvar^a$ are treated as independent spatial dependent variables.
Vanishing of the whole family characterizes a total spatial divergence locally;
the notation $\widehat{\sg E}_{\depvar^a}(f)=0$ means vanishing for every $k$.
The operators act componentwise on vector expressions.
These equations \eqref{eq:curl.conslaw.detsys} are necessary
(but not sufficient) for the multiplier functions to determine
conservation law identities of the form \eqref{eq.surf.conslaw.ident1}--\eqref{eq.surf.conslaw.ident2}. 

Each solution of the full compatibility identity \eqref{eq:curl.conslaw.deteq}
yields a local surface-flux conservation law through
Proposition~\ref{prop:operator.characteristics}.

\medskip\noindent\textbf{Example.}
For Maxwell's equations \eqref{eq:Max:vac} in Cartesian coordinates,
written as the system
$\vec{G}_E = \vec{E}_t - c\,\curl \vec{B}$,
$\vec{G}_B = \vec{B}_t + c\,\curl \vec{E}$,
$G^7 = \div\vec{E}$,
$G^8 = \div\vec{B}$
for $\depvars = (\vec{E},\vec{B})$,
take $\xi=0$ and consider multipliers depending on $t$ and $\vec{x}$ only.
Writing the multiplier combinations as
$Q_\sigma G^\sigma = \vec{q}_E\cdot\vec{G}_E + \vec{q}_B\cdot\vec{G}_B + q_7 G^7 + q_8 G^8$
and
$\vec{Q}_\sigma G^\sigma
= \underline{M}_E\, \vec{G}_E + \underline{M}_B\, \vec{G}_B + \vec{m}_7\, G^7 + \vec{m}_8\, G^8$
with matrices $\underline{M}_E$, $\underline{M}_B$,
the determining equation \eqref{eq:curl.conslaw.deteq}
splits with respect to the independent derivatives of $\vec{E}$, $\vec{B}$:
the coefficients of $\vec{E}_{tt}$, $\vec{B}_{tt}$ force
$\vec{q}_E=\vec{q}_B=0$;
those of the mixed derivatives $\partial_{x^i} E^j_t$, $\partial_{x^i} B^j_t$ force
$\underline{M}_E=q_7\,\idtens$, $\underline{M}_B=q_8\,\idtens$;
those of $\grad(\div\vec{E})$, $\grad(\div\vec{B})$ force
$\vec{m}_7=\vec{m}_8=0$,
the curl blocks dropping out through $\div\curl\equiv0$;
and the remaining first-order terms force
$q_7$, $q_8$ to be constants.
Reconstructed from the identities
\eqref{eq.surf.conslaw.ident1}--\eqref{eq.surf.conslaw.ident2},
the two resulting conservation laws are precisely
the electric and magnetic surface-flux conservation laws
$(\vec{\ConsLawDensity},\vec{\ConsLawFlux}) = (\vec{E},\,-c\vec{B})$
and $(\vec{B},\,c\vec{E})$
of Section~\ref{sec:surf:staticdomain}.
Moreover, within this multiplier ansatz,
Maxwell's equations admit no other independent nontrivial surface-flux conservation laws.

\subsubsection{Local gradient-type conservation laws}

For static gradient-type laws, take the special ansatz
$\vec q=\vec\xi(\vec x)+\sum_\sigma\vec Q_\sigma G^\sigma$
and $\vec b=\sum_\sigma\vec P_\sigma G^\sigma$ in
Proposition~\ref{prop:operator.characteristics}.
Here $\vec\xi$ is prescribed and field-independent, with $\div\vec\xi=0$,
and both multiplier families $\vec Q_\sigma$ and $\vec P_\sigma$ are vector-valued.
After the determining identities below are imposed, the reconstructed identities are
\beq\label{eq.curv.conslaw.ident1}
\curl\vec{\ConsLawDensity}
 =\vec\xi(\vec x)+\vec Q_1G^1+\cdots+\vec Q_NG^N
\eeq
and
\beq\label{eq.curv.conslaw.ident2}
\partial_t\vec{\ConsLawDensity}+\grad\ConsLawFlux
 =\vec P_1G^1+\cdots+\vec P_NG^N.
\eeq
This is again a sufficient algebraic ansatz, not a universal reduction
of the finite-operator representation.
Its complete determining identities are
\beq\label{eq:grad.conslaw.closure}
\div\big(\vec\xi+\vec Q_1G^1+\cdots+\vec Q_NG^N\big)\equiv0
\eeq
and
\beq\label{eq:grad.conslaw.deteq}
\partial_t(\vec Q_1G^1+\cdots+\vec Q_NG^N)
 -\curl(\vec P_1G^1+\cdots+\vec P_NG^N)\equiv0.
\eeq
They imply the necessary Euler tests
\beq\label{eq:grad.conslaw.detsys}
\begin{aligned}
&\widehat{\sg E}_{\depvar^a}\!\left(J\sum_\sigma Q_\sigma^{\,i}G^\sigma\right)\equiv0,\\
&\sg E_{\depvar^a}\!\left(g_{ij}\sum_\sigma P_\sigma^{\,j}G^\sigma\right)\equiv0,
\qquad a=1,\ldots,m,\quad i=1,2,3,
\end{aligned}
\eeq
where the spatial family is defined by \eqref{eq:spatial.euler.family}.
The first test uses the contravariant curl components,
$J(\curl\vec{\ConsLawDensity})^i=\epsilon^{ijk}D_j\ConsLawDensity_k$.
The second uses covariant components because
$g_{ij}(D_t\vec{\ConsLawDensity}+\grad\ConsLawFlux)^j
=D_t\ConsLawDensity_i+D_i\ConsLawFlux$.
Thus the gradient-law test requires lowering the vector index, rather than
multiplying every vector component by $J$.
Every solution of the full identities \eqref{eq:grad.conslaw.closure}--\eqref{eq:grad.conslaw.deteq}
locally reconstructs a circulatory conservation law; the Euler tests alone are not sufficient.

\medskip\noindent\textbf{Example.}
For Maxwell's equations in Cartesian coordinates, with $\vec\xi=0$ and multipliers depending on $t$ and $\vec{x}$ only,
the determining equation \eqref{eq:grad.conslaw.deteq}
splits analogously and admits only the zero solution.
Thus, Maxwell's equations possess no nontrivial circulatory conservation laws
within this algebraic multiplier ansatz. 

\medskip
Derivative-dependent algebraic multipliers and higher-order operator
characteristics can be sought using the full determining identities.
Their classification in particular PDE systems is left for future work.

\subsubsection{Multipliers for moving conservation laws}

For volumetric laws the relation \eqref{eq:CLflux:mov.vol}
allows the ordinary multiplier equations to be used without change.
The static surface and curve equations likewise apply to the simplifying
subclasses \eqref{eq:surfdens.divfree} and \eqref{eq:curvdens.curlfree},
using the flux relations \eqref{eq:CLflux:mov.surf} and \eqref{eq:CLflux:mov.curv}.

For the general moving laws, the determining identities are instead
\eqref{eq:operator.surface} and \eqref{eq:operator.curve}.
Their additional transport terms retain the divergence or curl of the density,
which need not vanish or be prescribed independently of the solution.
As a check, the electric Maxwell law has $q=G^7=\div\vec E$ and
$\vec b=\vec G_E+\VelVec G^7$. Then
\beq
D_tG^7+\div(\VelVec G^7)-\div(\vec G_E+\VelVec G^7)
 =D_tG^7-\div\vec G_E\equiv0.
\eeq
Reconstruction gives $\vec{\ConsLawDensity}=\vec E$ and
$\widetilde{\vec{\ConsLawFlux}}=-c\vec B-\vec E\times\VelVec$,
in agreement with \eqref{eq:Max:CL:moving:Ampere}.

\section{Discussion}\label{sec:conclusion}

Reynolds' transport theorem in fluid mechanics provides a generalization of
divergence-type conservation laws to material moving volumes
in three spatial dimensions. 
A main result of the present work has been to give the counterpart transport formulas
in both local and global forms for curl-type and gradient-type conservation laws
on moving surfaces and moving curves, respectively.
These transport formulae have been stated in the familiar vector calculus notation
used in fluid mechanics and applied mathematics.
In addition, the derivations have used an arbitrary velocity field,
which is not necessarily required to be a dynamical variable such as the fluid velocity 
in an underlying physical PDE system. 
The resulting transport formulae have also been shown to agree
with the more mathematical Lie-derivative formulation of transport theorems
\cite{Anc.2013,Niv.Cor.Kai.Sch.Nao} 
for moving submanifolds in a $n$-dimensional manifold specialized to $\mathbb{R}^3$. 

Analogous results apply in two spatial dimensions,
where plane domains play the role of both volumes and surfaces,
giving divergence-type conservation laws on plane domains
and gradient-type conservation laws on plane curves,
with the curl operator reducing to a dual type of divergence operator. 
The moving-domain generalizations given in Section~\ref{sec:CL:moving}
carry over directly. 
A notable two-dimensional illustration is provided by
inviscid incompressible fluid flow in a flat plane,
where the scalar vorticity $\omega = \p_x u^2 - \p_y u^1$
is materially conserved, $\omega_t + \vec{u}\cdot\nabla\omega = 0$.
As a consequence, for an arbitrary function $f$,
\beq
\p_t f(\omega) + \div\big( f(\omega)\, \vec{u}\big) = 0
\eeq
is a moving planemetric conservation law with zero moving flux,
using the fluid velocity as the transport velocity field. 
The resulting global conserved quantities
$\int_{\V(t)} f(\omega)\, dA$
on material plane domains include the enstrophy ($f(\omega)=\omega^2$)
and all higher moments of the vorticity.
A full discussion of moving conservation laws in two spatial dimensions
will be given in forthcoming work. 

As explained in Section~\ref{sec:multipliers},
all of the different types of local conservation laws on static and moving domains
for the PDE systems considered there can be characterized by multipliers
and their differential-operator extensions.
The multiplier method for divergence-type conservation laws is well known
(see e.g. Ref.~\cite{olver1993applications,bluman2010applications,Anc-review}),
and a main result of the present work has been to extend it
to curl-type and gradient-type conservation laws.
This extension provides a systematic method for finding local conservation laws
of all three types.
For a prescribed transport field, the determining equations are linear
in the unknown multipliers and differential-operator coefficients.
Computations can be organized by specifying their differential order
and dependence on the field variables and their derivatives.
In simple cases, calculations may be done by hand;
in more complex situations, symbolic computation becomes essential.
For divergence-type conservation laws,
the multiplier method has been implemented in various computer algebra systems
(see, e.g., \cite{cheviakov2007gem, cheviakov2010computation, cheviakov2010symbolic} and references therein).
A corresponding implementation of the extended method for curl-type and gradient-type conservation laws will be carried out in future work. 
The multiplier framework has also recently been extended to
\emph{approximate} conservation laws of PDE systems involving a small parameter \cite{Che.Tar.Yang.2025};
the development of counterparts
for approximate curl-type and gradient-type conservation laws 
will be pursued elsewhere.

\section*{Acknowledgements}

A.F.S. and S.C.A. thank NSERC for research support through individual Discovery grants.

\begin{appendix}

\section{Leibniz rule for moving domains}\label{sec:proof}

A rigorous proof of the transport formula \eqref{eq:Leibniz.vol},
which is commonly referred to as the \emph{Leibniz rule} on moving volumes,
will now be given.

Let $\mathcal{V}(t)$ be a material region 
which is carried by a velocity field $\VelVec(\vec{x},t)$.
In Cartesian coordinates, use an orientation-preserving parametrization 
$\vec{x} = \vec{X}(a, b, c, t)$
in terms of parameters 
$(a, b,c) \in U \subset \mathbb{R}^3$,
where the parameter domain $U$ is fixed, 
with $\partial_t\vec{X} = \VelVec(\vec{X}, t)$.
The pulled-back volume element on $U$ is $\vec{X}_a\cdot (\vec{X}_b\times\vec{X}_c)\,da\,db\,dc$.
Taking the time derivative of $\vec{X}_a\cdot (\vec{X}_b\times\vec{X}_c) = \epsilon_{ijk} X^i_a X^j_b X^k_c$ yields
\beq
\begin{aligned}
\epsilon_{ijk} \partial_t(X^i_a X^j_b X^k_c) 
& = 3\epsilon_{ijk} v^i(\vec{X},t)_{[a}^{\mathstrut} X^j_{b} X^k_{c]} 
= 3\epsilon_{ijk} \partial_{x^l}v^i\, X^l_{[a} X^j_{b} X^k_{c]} 
= \partial_{x^l}v^l \epsilon_{ijk} X^i_{[a} X^j_{b} X^k_{c]} \\
& = \div\VelVec\; (\vec{X}_a\cdot (\vec{X}_b\times\vec{X}_c)) . 
\end{aligned}
\eeq
Therefore, the time derivative of the volume integral of a density
given by a differential scalar function $\ConsLawDensity$ is given by 
\beq
\begin{aligned}
  \frac{d}{dt} \mathop{\int}_{\mathcal{V}(t)} \ConsLawDensity  \;dV
  & =
  \frac{d}{dt} \mathop{\int}_{U} \ConsLawDensity\; (\vec{X}(a,b,c,t),t) \vec{X}_a\cdot (\vec{X}_b\times\vec{X}_c)\,da\,db\,dc 
  \\[2ex]
  & = 
  \mathop{\int}_{U} (\partial_t \ConsLawDensity + \VelVec\cdot\grad\ConsLawDensity +\ConsLawDensity \div\VelVec)\; (\vec{X}_a\cdot (\vec{X}_b\times\vec{X}_c)) \,da\,db\,dc 
  \\[2ex]
  & = \mathop{\int}_{\mathcal{V}(t)} \left( \partial_t\ConsLawDensity  + \div (\ConsLawDensity\, \VelVec) \right)\;dV,
\end{aligned}
\eeq
where the integrands over $U$ are evaluated at $\vec{x}=\vec{X}(a,b,c,t)$,
and the product rule and change of variables have been used in the last step.
This completes the derivation of formula \eqref{eq:Leibniz.vol}.

The Leibniz rule on a moving surface,
described by the transport formula \eqref{eq:Leibniz.surf.gen},
follows by differentiating a transported parametrization.
In Cartesian coordinates, write a surface patch as
$\vec{x}=\vec{X}(a,b,t)$, where $a,b$ are surface parameters,
$(a,b)\in U\subset\mathbb{R}^2$, and the parameter domain $U$ is fixed,
with $\partial_t\vec{X}=\VelVec(\vec{X},t)$.
The oriented area vector is
$\vec{N}=\vec{X}_a\times\vec{X}_b$, so that
$\dS=\vec{N}\,da\,db$.
Differentiating the two tangent vectors gives
\beq
\partial_t N^i=(\partial_{x^j}v^j)N^i-(\partial_{x^i}v^j)N^j,
\eeq
where the velocity derivatives are evaluated at $\vec{X}(a,b,t)$.
Consequently,
\beq
\begin{aligned}
\dfrac{d}{dt}\int_{\S(t)}\vec{\ConsLawDensity}\cdot\dS
&=\dfrac{d}{dt}\int_U
\vec{\ConsLawDensity}(\vec{X}(a,b,t),t)\cdot\vec{N}\,da\,db\\
&=\int_{\S(t)}\Big(
\partial_t\vec{\ConsLawDensity}
+(\VelVec\cdot\nabla)\vec{\ConsLawDensity}
+(\div\VelVec)\vec{\ConsLawDensity}
-(\vec{\ConsLawDensity}\cdot\nabla)\VelVec
\Big)\cdot\dS.
\end{aligned}
\eeq
The vector identity
\beq
\curl(\vec{\ConsLawDensity}\times\VelVec)
=(\VelVec\cdot\nabla)\vec{\ConsLawDensity}
-(\vec{\ConsLawDensity}\cdot\nabla)\VelVec
+(\div\VelVec)\vec{\ConsLawDensity}
-(\div\vec{\ConsLawDensity})\VelVec
\eeq
therefore yields
\beq\label{eq:Leibnitz.surf}
\begin{aligned}
\dfrac{d}{dt}\int_{\S(t)}\vec{\ConsLawDensity}\cdot\dS
&=\int_{\S(t)}\Big(
\partial_t\vec{\ConsLawDensity}
+(\div\vec{\ConsLawDensity})\VelVec
+\curl(\vec{\ConsLawDensity}\times\VelVec)
\Big)\cdot\dS\\
&=\int_{\S(t)}\Big(
\partial_t\vec{\ConsLawDensity}
+(\div\vec{\ConsLawDensity})\VelVec
\Big)\cdot\dS
+\oint_{\p\S(t)}(\vec{\ConsLawDensity}\times\VelVec)\cdot\dl,
\end{aligned}
\eeq
where Stokes' theorem gives the second equality.
For a surface covered by several patches, the internal boundary integrals cancel.
This proves \eqref{eq:Leibniz.surf.gen};
condition \eqref{eq:surfdens.divfree} gives its reduction
\eqref{eq:Leibniz.surf}.

For a moving curve, use a transported parametrization
$\vec{x}=\vec{X}(a,t)$, $a_0\leq a\leq a_1$,
with $\partial_t\vec{X}=\VelVec(\vec{X},t)$ and
$\dl=\vec{X}_a\,da$.
Differentiating the line integral gives
\beq
\begin{aligned}
\dfrac{d}{dt}\int_{\C(t)}\vec{\ConsLawDensity}\cdot\dl
&=\int_{a_0}^{a_1}\Big(
\partial_t\ConsLawDensity^i
+v^j\partial_{x^j}\ConsLawDensity^i
+\ConsLawDensity^j\partial_{x^i}v^j
\Big)X_a^i\,da\\
&=\int_{\C(t)}\Big(
\partial_t\vec{\ConsLawDensity}
+(\curl\vec{\ConsLawDensity})\times\VelVec
+\grad(\vec{\ConsLawDensity}\cdot\VelVec)
\Big)\cdot\dl.
\end{aligned}
\eeq
The fundamental theorem for line integrals now yields
\beq\label{eq:Leibnitz.curv}
\dfrac{d}{dt}\int_{\C(t)}\vec{\ConsLawDensity}\cdot\dl
=\int_{\C(t)}\Big(
\partial_t\vec{\ConsLawDensity}
+(\curl\vec{\ConsLawDensity})\times\VelVec
\Big)\cdot\dl
+\big(\vec{\ConsLawDensity}\cdot\VelVec\big)
\Big|_{\vec{x}_a(t)}^{\vec{x}_b(t)}.
\eeq
This is \eqref{eq:Leibniz.curv.gen}.
The endpoint term involves the full endpoint velocity;
it vanishes for a closed curve.
Condition \eqref{eq:curvdens.curlfree} gives the reduction
\eqref{eq:Leibniz.curv}.

\section{Leibniz rule via Lie derivatives}\label{sec:liederivative}

Here $\nabla$ denotes the Levi-Civita covariant derivative of the Euclidean
metric, expressed in the chosen coordinates.
We use the determinant normalization of the wedge product
\cite[Eq.~(1.49), p.~55]{olver1993applications}; thus a differential $r$-form
with components $\omega_{i_1\cdots i_r}$ is written as
$\omega=\frac{1}{r!}\omega_{i_1\cdots i_r}\,dx^{i_1}\wedge\cdots\wedge dx^{i_r}$,
with summation over all repeated indices.

There is an alternative way to prove the Leibniz rule for moving domains,
which uses Lie derivatives and differential forms as shown in Ref.~\cite{Anc.2013,Niv.Cor.Kai.Sch.Nao}. 
This proof will now be outlined in the index notation summarized in Appendix~\ref{sec:curvilin.coords},
which will establish its equivalence with the vector-calculus formulas 
\eqref{eq:Leibniz.vol}, \eqref{eq:Leibniz.surf.gen}, \eqref{eq:Leibniz.curv.gen}
for moving volumes, surfaces, and curves.
Note that these formulae are coordinate invariant. 

To begin,
the volume integral of density $\ConsLawDensity$ on a moving domain $\mathcal{V}(t)$
can be expressed as
\beq
\mathop{\int}_{\mathcal{V}(t)} \ConsLawDensity\, \tfrac{1}{3!}\varepsilon_{jkl}\,dx^j\wedge dx^k \wedge dx^l
\eeq
where $\varepsilon_{jkl}$ denotes the components of the volume form
in an arbitrary coordinate system $(x^1,x^2,x^3)$. 
Using the Lie derivative,
the time derivative of this integral is given by
\beq\label{integral.vol.lieder}
\dfrac{d}{dt} \mathop{\int}_{\mathcal{V}(t)} \ConsLawDensity\, \tfrac{1}{3!}\varepsilon_{jkl}\, dx^j\wedge dx^k \wedge dx^l
= \mathop{\int}_{\mathcal{V}(t)} \Dop_t (\ConsLawDensity\, \tfrac{1}{3!}\varepsilon_{jkl}\,dx^j\wedge dx^k \wedge dx^l)
\eeq
where
\beq
\Dop_t = \partial_t + \lieder{\VelVec}
\eeq
denotes the material Lie derivative.
It acts as the ordinary material derivative on differential scalar functions,
\beq
\Dop_t\ConsLawDensity = \partial_t\ConsLawDensity + \VelVec\cdot\nabla \ConsLawDensity . 
\eeq
On the volume form,
the material Lie derivative is given by
\beq
\Dop_t (\tfrac{1}{3!}\varepsilon_{jkl}\,dx^j\wedge dx^k \wedge dx^l)
= \lieder{\VelVec}(\tfrac{1}{3!}\varepsilon_{jkl}\,dx^j\wedge dx^k \wedge dx^l)
= (\div\VelVec) \tfrac{1}{3!}\varepsilon_{jkl}\,dx^j\wedge dx^k \wedge dx^l . 
\eeq
Hence, one has
\beq\begin{aligned}
\Dop_t (\ConsLawDensity\, \tfrac{1}{3!}\varepsilon_{jkl}\,dx^j\wedge dx^k \wedge dx^l)
& = \big( \partial_t\ConsLawDensity + \VelVec\cdot\nabla \ConsLawDensity
+ (\div\VelVec)\ConsLawDensity \big) \tfrac{1}{3!}\varepsilon_{jkl}\,dx^j\wedge dx^k \wedge dx^l
\\& 
= \big( \partial_t\ConsLawDensity + \div(\ConsLawDensity\VelVec) \big) \tfrac{1}{3!}\varepsilon_{jkl}\,dx^j\wedge dx^k \wedge dx^l, 
\end{aligned}\eeq
and thus the time derivative integral \eqref{integral.vol.lieder}
satisfies the formula \eqref{eq:Leibniz.vol}.

Next, the surface integral of a vector density $\ConsLawDensity^l$
on a moving surface $\mathcal{S}(t)$
can be expressed as
\beq
\mathop{\int}_{\mathcal{S}(t)} \ConsLawDensity^i\, \hat{n}_i\, \hat{n}^l\, \tfrac{1}{2}\varepsilon_{jkl}\, dx^j \wedge dx^k
= \mathop{\int}_{\mathcal{S}(t)} \ConsLawDensity^l\, \tfrac{1}{2}\varepsilon_{jkl}\, dx^j \wedge dx^k
\eeq
where $\hat{n}^i$ is the unit normal vector of the surface.
The pullback of this 2-form to the surface yields a projection onto the normal component of the vector density. 
The time derivative of the surface integral is then given by
\beq\label{integral.surf.lieder}
\dfrac{d}{dt} \mathop{\int}_{\mathcal{S}(t)} \ConsLawDensity^l\, \tfrac{1}{2}\varepsilon_{jkl}\, dx^j \wedge dx^k
= \mathop{\int}_{\mathcal{S}(t)} \Dop_t (\ConsLawDensity^l\, \tfrac{1}{2}\varepsilon_{jkl} \, dx^j \wedge dx^k)
\eeq
where the material Lie derivative acts on the vector function
$\ConsLawDensity^l\dfrac{\partial}{\partial x^l}$ by 
\beq
\Dop_t \Big(\ConsLawDensity^l \dfrac{\partial}{\partial x^l} \Big)
= \partial_t \ConsLawDensity^l \dfrac{\partial}{\partial x^l}
+ \lieder{\VelVec}\Big( \ConsLawDensity^l \dfrac{\partial}{\partial x^l} \Big)
= \partial_t \ConsLawDensity^l \dfrac{\partial}{\partial x^l}
+ \Big( \VelVec\cdot\nabla\ConsLawDensity^l - \vec{\ConsLawDensity}\cdot\nabla v^l \Big)\dfrac{\partial}{\partial x^l} . 
\eeq
Hence, one has
\beq
\Dop_t (\ConsLawDensity^l\, \tfrac{1}{2}\varepsilon_{jkl} \, dx^j \wedge dx^k)
= 
\big( \partial_t \ConsLawDensity^l + \VelVec\cdot\nabla\ConsLawDensity^l
- \vec{\ConsLawDensity}\cdot\nabla v^l + (\div\VelVec) \ConsLawDensity^l \big) \tfrac{1}{2}\varepsilon_{jkl}\, dx^j \wedge dx^k . 
\eeq
Next use the curl identity 
\beq
\curl(\vec{\ConsLawDensity}\times \VelVec)^l
= \VelVec\cdot\nabla\ConsLawDensity^l - \vec{\ConsLawDensity}\cdot\nabla v^l
+ (\div\VelVec) \ConsLawDensity^l  - (\div\vec{\ConsLawDensity}) v^l , 
\eeq
which gives
\beq
\Dop_t (\ConsLawDensity^l\, \tfrac{1}{2}\varepsilon_{jkl}\, dx^j \wedge dx^k)
= \big( \partial_t \ConsLawDensity^l + \varepsilon^l{}_{mn} \nabla^m (\vec{\ConsLawDensity}\times \VelVec)^n + (\div\vec{\ConsLawDensity}) v^l \big) \tfrac{1}{2}\varepsilon_{jkl}\, dx^j \wedge dx^k . 
\eeq
Thus, the time derivative integral \eqref{integral.surf.lieder}
satisfies the formula \eqref{eq:Leibniz.surf.gen}.

Last, the line integral of a vector density $\ConsLawDensity^l$
on a moving curve $\mathcal{C}(t)$ 
can be expressed as
\beq
\mathop{\int}_{\mathcal{C}(t)} \ConsLawDensity^i\, \hat{\ell}_i\, \hat{\ell}_j\, dx^j
= \mathop{\int}_{\mathcal{C}(t)} \ConsLawDensity_j \, dx^j
\eeq
where $\hat{\ell}^i$ is the unit tangent vector of the curve,
and $\ConsLawDensity_j\, dx^j= g_{jk} \ConsLawDensity^k\,dx^j$ is the 1-form
corresponding to the vector density. 
The pullback of this 1-form to the curve in the line integral
applies a projection 
$\ConsLawDensity_k \hat{\ell}_j\, \hat{\ell}^k = \hat{\ell}_j \ConsLawDensity^k\hat{\ell}_k$. 
Consequently, the time derivative of the line integral is given by
\beq\label{integral.curv.lieder}
\dfrac{d}{dt} \mathop{\int}_{\mathcal{C}(t)} \ConsLawDensity_j \, dx^j
= \mathop{\int}_{\mathcal{C}(t)} \Dop_t (\ConsLawDensity_j \, dx^j)
\eeq
where the action of the material Lie derivative consists of 
\beq
\Dop_t (\ConsLawDensity_j \, dx^j)
= \partial_t (\ConsLawDensity_j\, dx^j) + \lieder{\VelVec} (\ConsLawDensity_j\, dx^j)
= (\partial_t \ConsLawDensity_j +\VelVec\cdot\nabla\ConsLawDensity_j 
+ (\nabla_j \VelVec)\cdot \vec{\ConsLawDensity} )\, dx^j. 
\eeq
Then the curl identity
\beq
(\VelVec\times \curl \vec{\ConsLawDensity})^l
= (\nabla^l \vec{\ConsLawDensity})\cdot \VelVec  - \VelVec\cdot\nabla \ConsLawDensity^l
\eeq
gives
\beq
\begin{aligned}
\Dop_t (\ConsLawDensity_j \, dx^j)
&= \big(\partial_t \ConsLawDensity_j - (\VelVec\times \curl \vec{\ConsLawDensity})_j
+ (\nabla_j \vec{\ConsLawDensity})\cdot \VelVec
+ (\nabla_j \VelVec)\cdot \vec{\ConsLawDensity}\big)\, dx^j\\
&= \big(\partial_t \ConsLawDensity_j + ((\curl \vec{\ConsLawDensity})\times \VelVec)_j
+ \nabla_j(\vec{\ConsLawDensity}\cdot\VelVec)\big) \, dx^j.
\end{aligned} 
\eeq
Thus, the time derivative integral \eqref{integral.curv.lieder}
satisfies the formula \eqref{eq:Leibniz.curv.gen}.

\section{Curvilinear coordinates}\label{sec:curvilin.coords}

Consider an arbitrary positively oriented, time-independent curvilinear coordinate system $\vec{x}=(x^1,x^2,x^3)$
in $\mathbb{R}^3$.
It will be convenient to use index notation $\vec{x} = x^i$, $i=1,2,3$,
and employ the summation convention that any repeated index is summed over.

Let $a^i$ and $b^i$ denote the curvilinear components of vectors $\vec{a}$ and $\vec{b}$.
The dot product and cross product of vectors 
are given by
\beq
\vec{a}\cdot\vec{b} = g_{ij} a^i b^j 
\eeq
and
\beq
(\vec{a}\times\vec{b})^i = g^{ij}\varepsilon_{jkl} a^k b^l 
\eeq
where $g_{ij}$ and $\varepsilon_{jkl}$ denote the components of
the metric tensor and volume form in the given coordinate system.
Note that, in Cartesian coordinates, 
$g_{ij} = \delta_{ij}$ is the Kronecker symbol
and $\varepsilon_{jkl} = \epsilon_{jkl}$ is the Levi-Civita symbol.
In any other coordinates, these tensors are given by
the covariant transformation formula in terms of the Jacobian
$J= \frac{\partial(x,y,z)}{\partial(x^1,x^2,x^3)} = (J_i^1,J_i^2,J_i^3)$:
\beq
g_{ij} = J_i{}^1 J_j{}^1 + J_i{}^2 J_j{}^2 + J_i{}^3 J_j{}^3
\eeq
with $\delta_{11}=\delta_{22}=\delta_{33}=1$
being the non-zero components of the Kronecker symbol;
\beq
\varepsilon_{ijk} = 6 J_{[i}{}^1 J_j{}^2 J_{k]}{}^3  = \sqrt{\det(g)}\,\epsilon_{ijk} 
\eeq 
where the square brackets denote total antisymmetrization of the indices,
with $\epsilon_{123}=\epsilon_{231}=\epsilon_{312}=-\epsilon_{213}=-\epsilon_{132}=-\epsilon_{321}=1$
being the non-zero components of Levi-Civita symbol. 

The divergence and curl of a vector field $\vec{f}(\vec{x})$ are given by
\beq
\div\vec{f} = \frac{1}{\sqrt{\det(g)}} \partial_{x^i} (\sqrt{\det(g)} f^i) 
\eeq
and
\beq
(\curl\vec{f})^i  = \frac{1}{\sqrt{\det(g)}} \epsilon^{ijk} \partial_{x^j} ( g_{kl} f^l )
\eeq
with $\epsilon^{123}=\epsilon_{123}$ and so on. 
The gradient of a scalar field (or function) $f(\vec{x})$ is given by
\beq
(\grad f)^i  = g^{ij} \partial_{x^j} f
\eeq
where $g^{ij}$ is inverse of $g_{ij}$. 

The Laplacian of a scalar field is given by
\beq
\Delta f = \div\grad f = \frac{1}{\sqrt{\det(g)}} \partial_{x^i} (\sqrt{\det(g)} g^{ij} \partial_{x^j} f)
\eeq
while the Laplacian of a vector field $\vec{f}(\vec{x})$ is given by 
\beq\label{eq:vector:Lapl}
\Delta \vec{f} = \grad\, ( \div \vec{f}) - \curl\, (\curl \vec{f})
\eeq
whose explicit component form follows by composing
the divergence, gradient, and curl formulas above.

\end{appendix}

\bibliographystyle{ieeetr}
\bibliography{Mov3DCLs-AS-SA-v14}

\end{document}